\documentclass{amsart}
\usepackage[a4paper, total={6in, 8in}]{geometry} 	
\usepackage[utf8]{inputenc}
\usepackage[english]{babel}
\usepackage{mypreamble}

\title{The Goresky-Hingston Coproduct on Based Loop Spaces}
\author{Jonathan Clivio}
\begin{document}
\begin{abstract}
    We construct a lift of the Goresky-Hingston coproduct on the based loop space. Using this lift, we produce formulas for the interaction of the coproduct with the $\pi_1$-action and with maps of manifolds. These results lets us compute the coproduct on the based loop space for manifolds of the form $S^n/G$.
\end{abstract}

\maketitle

\section{Introduction}

In this article, we study the Goresky-Hingston coproduct on the based loop space of a pointed, closed, oriented, $n$-dimensional manifold $(M,x_0)$:
\begin{align*}
    \lor_{\mathrm{GH}}\colon H_*(\Omega M,x_0)\to H_{1-n+*}((\Omega M,x_0)\times (\Omega M,x_0)).
\end{align*}
This coproduct was defined in \cite{goresky2009loop} and \cite{sullivan2004open}.

In this article, we construct a slightly different coproduct, the \textit{avoiding-stick coproduct}
\begin{align*}
    \lor_{\Omega M}\colon H_*(\Omega M)\to H_{1-n+*}(\Omega M\times \Omega M).
\end{align*}
The avoiding-stick coproduct is a lift of the Goresky-Hingston coproduct (see Proposition \ref{prop:a-s-coprod lifts GH}).

Because this coproduct does not have any relative terms, it is easier to work with. We prove an array of results for the avoiding-stick coproduct that directly descend to results for the Goresky-Hingston coproduct. The first result is an invariance result under a $\pi_1$-action:

\begin{IntroThm}[Corollary \ref{cor:pi1 invariance}]\label{introthm:pi1 invariance}
    Assume that the dimension is $n>1$. The avoiding-stick coproduct $\lor$ factors through the fixed points $H_*(\Omega M\times \Omega M)^{\pi_1(M,x_0)}$ where $\pi_1(M,x_0)$ acts on $\Omega M\times \Omega M$ via
    \begin{align*}
        [\tau]\cdot(\gamma,\gamma'):=(\gamma \cdot\tau, \tau^{-1} \cdot\gamma').
    \end{align*}
\end{IntroThm}
This immediately implies that the avoiding-stick and thus the Goresky-Hingston coproduct both are zero if $\pi_1(M,x_0)$ is infinite (see Corollary \ref{cor:infty pi1}).

Our second main result is a formula for how the avoiding-stick coproduct interacts with maps between manifolds $f\colon M\to N$. It generalizes a result in a note by Nancy Hingston \cite[Lemma 3]{hingston2010loop} to the case where $f$ is of any degree and $N$ has arbitrary fundamental group:
\begin{IntroThm}[Corollary \ref{cor:formula for f}]\label{introthm:formula for f}
    Let $f\colon (M,x_0)\to (N,y_0)$ be a map of closed, oriented, pointed manifolds of dimension $n>1$. Assume that $y_0$ is a regular value with $f^{-1}(y_0)=\{x_0,\dots, x_k\}$. Let $\tau_i$ be paths from $x_i$ to $x_0$. For the avoiding-stick coproducts $\lor_{\Omega M}$ and $\lor_{\Omega N}$, it holds
    \begin{align*}
        \lor_{\Omega N} \circ f_*=\sum_{i=0}^k \deg(f)_{x_i} (r_{[f_*(\tau_i)]}\times l_{[f_*(\tau_i)]^{-1}})\circ (f\times f)_*\circ \lor_{\Omega M}
    \end{align*}
    where $r_{\gamma}$ and $l_\gamma$ denote right and left multiplication with an element $\gamma\in \pi_1(N,y_0)$.
\end{IntroThm}

This has as immediate consequence that the avoiding-stick and Goresky-Hingston coproduct on the based loop space commute with degree $1$ maps (see Corollary \ref{cor:deg 1}). This is in contrast to the Goresky-Hingston coproduct on the free loop space (see \cite{naef2021string,naef2024simple,wahl2019invariance,kenigsberg2024obstructions}).

Applying Theorem \ref{introthm:formula for f} to the case where $f$ is a universal covering lets us describe the coproduct fully in terms of its universal covering:
\begin{IntroThm}[Remark \ref{rem:Omega X univ cov} and Proposition \ref{prop:GH via univ cov}]
    Let $(M,x_0)$ be a pointed, closed, oriented manifold with finite fundamental group and universal covering $p\colon (\widetilde{M},y_0)\to (M,x_0)$. Any element $\alpha\in H_*(\Omega M)$ can be uniquely written as
    \begin{align*}
        \alpha=\sum_{g\in \pi_1(M,x_0)} g\cdot p_*(\beta_g)
    \end{align*}
    for some $\beta_g\in H_*(\Omega \widetilde M)$. If the avoiding-stick coproduct on $H_*(\Omega \widetilde{M})$ is computed by
    \begin{align*}
        \lor_{\Omega\widetilde{M}}(\beta_g)=\sum_i\gamma_{i,g}\times \delta_{i,g},
    \end{align*}
    then it holds
    \begin{align*}
        \lor_{\Omega M}(\alpha)=\sum_{g,h\in \pi_1(M)} \sum_i  g\cdot p_*(\gamma_{i,g})\cdot h\times h^{-1}\cdot p_*(\delta_{i,g}).
    \end{align*}
\end{IntroThm}

Because we can compute the coproduct on spheres (Example \ref{exmp:spheres}), we get a full description of the coproduct on $M=S^n/G$ where $G$ is a finite group acting freely on $S^n$:

\begin{IntroThm}[Proposition \ref{prop:Omega Sn/G} and Corollary \ref{cor:coprod Omega Sn/G}]\label{introthm:Omega Sn/G}
    Let $G$ be a finite group that acts on $S^n$ freely and orientation-preserving for some $n\geq 2$ and let $M=S^n/G$. The Pontryagin ring on $\Omega M$ is given by
    \begin{align*}
        H_*(\Omega M)\cong \Z[G][x]
    \end{align*}
    where $x$ is a central element of degree $n-1$.
    
    The avoiding-stick coproduct extends to a map $H_{n-1+*}(\Omega M)\to H_*(\Omega M)\otimes H_*(\Omega M)$ which is computed as
    \begin{align*}
        \lor(gx^k)=\sum_{i+j=k-1}\sum_{h\in G} gh^{-1}x^i\otimes hx^j.
    \end{align*}
\end{IntroThm}

This coproduct helps us understand the Leray-Serre spectral sequence associated to the fibration:
\begin{center}
    \begin{tikzcd}
        \Omega M \ar[r] & \Lambda M\ar[d, "\mathrm{ev}_0"] \\
        & M
    \end{tikzcd}
\end{center}
where $\Lambda M:=\Map(S^1,M)$ denotes the free loop space. Theorem \ref{introthm:Omega Sn/G} lets us show that all differentials on the $E^r$-page for $r\geq 2$ are zero and we find the following result:
\begin{IntroThm}[Corollary \ref{cor:free loop space}]
    Let $G$ be a non-trivial, finite group that acts on $S^n$ freely and orientation-preserving for some $n\geq 2$ and let $M=S^n/G$. The homology of the free loop space $\Lambda M$ has connected components $\Lambda_{[g]}M$ corresponding to conjugacy classes $[g]$ in $G$. The homology of $\Lambda_{[g]}M$ fits in the following short exact sequence:
    \begin{center}
        \begin{tikzcd}
            0 \ar[r] & H_{i-{n-1}}(S^n/C_G(g)) \ar[r]& H_{i+k(n-1)}(\Lambda_{[g]} M) \ar[r] & H_i(S^n/C_G(g)) \ar[r] &0
        \end{tikzcd}
    \end{center}
    for $k\geq 0$ and $1\leq i \leq n$. This short exact sequence splits for all $i\neq n-1$.
\end{IntroThm}

In \cite{clivio2025GHviaDGMorse}, we use this description to give a description of the coproduct on $\Lambda (S^n/G)$ with real coefficients where the coproduct is completely determined by the coproduct on $\Omega (S^n/G)$. This is in contrast to the coproduct with integer coefficients which \enquote{knows} about the simple homotopy of $S^n/G$: in \cite{naef2021string}, Naef computes parts of the coproduct on the lens spaces which is an example of $S^3/(\Z/7\Z)$ to give an example of a coproduct which does not commute with homotopy equivalences.

\subsection{The Avoiding-Stick Coproduct}
We briefly sketch the idea behind the avoiding-stick coproduct and why it does not produce any relative terms. In Subsection \ref{subsec:a-s def}, this is done in full detail.\\

We first recall a sketch of the definition of the Goresky-Hingston coproduct (see Subsection \ref{subsec:GH coprod} for more details). Roughly speaking, the coproduct on $\alpha\in C_*(\Omega M)$ is computed as follows: multiplying with the fundamental class of the interval $I=[0,1]$ gives a class $I\times \alpha\in C_{*+1}(I\times \Omega M,\partial I\times \Omega M)$. We then \enquote{intersect} this class with the subspace of self-intersections: all pairs $(t,\gamma)\in I\times \Omega M$ with $\gamma(t)=\gamma(0)$. Finally we apply the cutting map on the space of self-intersections: a pair $(\gamma,t)$ gets sent to the pair of loops $(t,\gamma)\mapsto(\gamma|_{[0,t]},\gamma|_{[t,1]})\in \Omega M$.

\begin{figure}[H]
    \centering
    \includegraphics[width=0.25\linewidth]{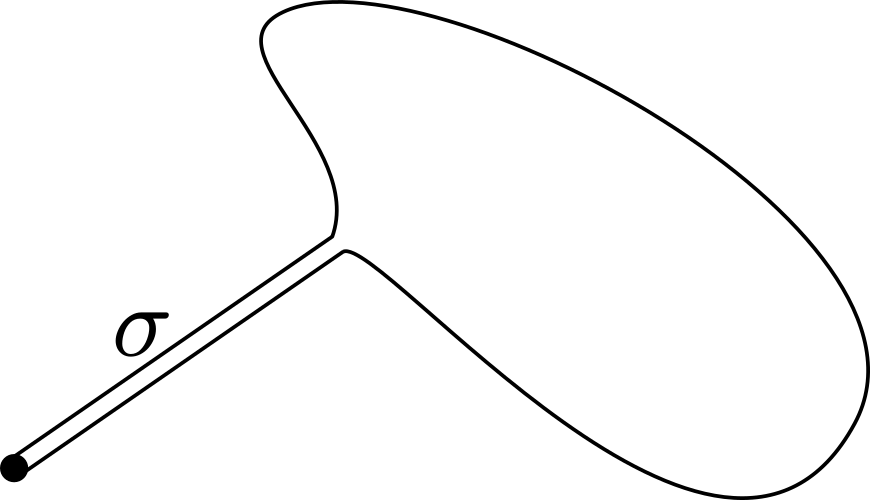}
    \caption{The space $\Omega^\sigma M$ consists of all loops that start and end with the path $\sigma$.}
    \label{fig:Omegasigma}
\end{figure}

All loops $\gamma\in \Omega M$ have trivial self-intersections for $t=0$ and $t=1$. Cutting at a trivial self-intersection gives a pair of loops $(\gamma_1,\gamma_2)\in \Omega M\times \Omega M$ where $\gamma_1$ or $\gamma_2$ is the constant loop $x_0$. This is why the Goresky-Hingston coproduct gives a map
\begin{align*}
    C_*(\Omega M)\to C_{1-n+*}(\Omega M\times \Omega M,\{x_0\}\times \Omega M\cup \Omega M\times \{x_0\}).
\end{align*}
The first advantage of the avoiding-stick coproduct is that it isolates any trivial self-intersection. It is defined analogously to the Goresky-Hingston coproduct but only on a subspace of $\Omega M$: we fix a path $\sigma\colon [0,\varepsilon]\to M$ a stick such that $\sigma(t)\neq x_0$ for $t>0$ and $\sigma(0)=x_0$. We consider the subspace $\Omega^\sigma M\subseteq \Omega M$ of all paths that follow $\sigma$ on $[0,\varepsilon]$ and $[1-\varepsilon,1]$ (see Figure \ref{fig:Omegasigma}). This subspace is homotopic to the whole space $\Omega M$. We define the avoiding-stick coproduct on $\Omega^\sigma M$ instead of $\Omega M$.
Working with the space $\Omega^\sigma M$ has the advantage that the trivial self-intersections are isolated because for all $t\in (0,\varepsilon]\cup [1-\varepsilon,1)$ and $\gamma\in \Omega^\sigma M$, we have $\gamma(t)\neq x_0$. All non-trivial self-intersections are at times $t\in (\varepsilon,1-\varepsilon)$. A cycle $\alpha\in C_*(I\times \Omega^\sigma M,\partial I\times \Omega^\sigma M)$ which represents some non-trivial self-intersections is thus also a non-relative cycle. Running the rest of the machinery of the Goresky-Hingston coproduct then produces no relative terms.

The avoiding-stick coproduct thus gives a map
\begin{align*}
    C_*(\Omega M)\cong C_*(\Omega^\sigma M)\to C_{1-n+*}(\Omega M\times \Omega M)
\end{align*}
which lifts the Goresky-Hingston coproduct.

Another advantage of the space $\Omega^\sigma M$ is that we have time at both ends of the loop to manipulate the loops. We can take advantage of this in the following way: let $H\colon [0,\varepsilon]\times M\to M$ be a $1$-parameter family of homeomorphisms such that $H(0,x)=x$ and $\widetilde{\sigma}(t)=H(t,\sigma(t))$ is also a stick avoiding $x_0$. We define the following homotopy equivalence:
\begin{align*}
    \Phi \colon \Omega^\sigma M&\to \Omega^{\widetilde{\sigma}}M,\\
    (\gamma\colon [0,1]\to M)&\mapsto t\mapsto\begin{cases}
        H(t,\gamma(t)) &\text{if } 0\leq t\leq \varepsilon;\\
        H(\varepsilon,\gamma(t)) &\text{if } \varepsilon\leq t\leq 1-\varepsilon;\\
        H(1-t,\gamma(t)) &\text{if } 1-\varepsilon\leq t\leq 1.
    \end{cases}
\end{align*}
\begin{figure}[H]
    \centering
    \includegraphics[width=0.3\linewidth]{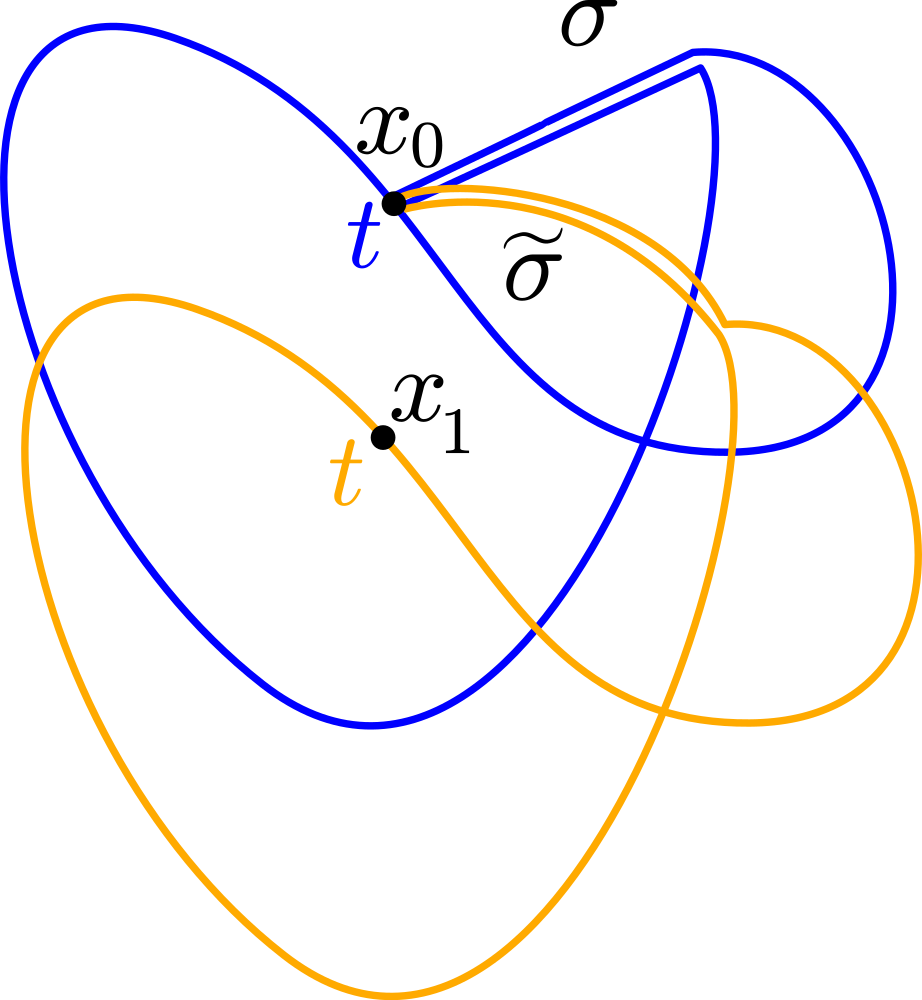}
    \caption{In this picture, the $1$-parameter family $H$ is given by translation downwards. The blue loop $\gamma$ is in $\Omega^\sigma M$ and has a self-intersection at $t$. Its image $\Phi(\gamma)$ in yellow then is in $\Omega^{\widetilde{\sigma}}M$ and at $t$ it crosses $x_1=\phi^\varepsilon(x_0)$.}
    \label{fig:hpty_Omegasigma}
\end{figure}
The map $\Phi$ sends the non-trivial self-intersections $(t,\gamma)\in (\varepsilon,1-\varepsilon)\times \Omega^\sigma M$ to all pairs $(t,\gamma')\in (\varepsilon,1-\varepsilon)\times \Omega^\sigma M$ with $\gamma'(t)=x_1:=H(\varepsilon,x_0)$ (see Figure \ref{fig:hpty_Omegasigma}). In fact, it gives a bijection. Thus we can detect self-intersections at the point $x_1$.  This is the insight that lets us prove Proposition \ref{prop:a-s commuting with homotopy} which is the key technical result unlocking the other statements of this article.

\addtocontents{toc}{\SkipTocEntry}
\subsection*{Organisation of the Paper}
In Section \ref{sec:various coproducts}, we define and compare the trivial coproduct, the Goresky-Hingston coproduct and the avoiding-sticks coproduct on the based loop space of a manifold. In Section \ref{sec:Cutting at Other Points}, we study what happens when we cut at a different point rather than the basepoint and relate it to the coproducts from the previous section. In Section \ref{sec:Maps and Coprod}, we describe how the different coproducts interact with maps of manifold. We apply this to the example of $S^n/G$ in Section \ref{sec:Sn/G}.

\addtocontents{toc}{\SkipTocEntry}
\subsection*{Conventions}
We denote by $I$ the closed interval $[0,1]$ and $\partial I=\{0,1\}$. For $\sigma\colon [a,b]\to X$ a path, we denote $\overline{\sigma}(t):=\sigma(b-t+a)$ for the path in reverse direction.

Let $M$ be a closed, oriented topological manifold of dimension $n$ and let $x_0\in M$ be a fixed basepoint. We denote by 
\begin{align*}
    \Omega M:=\{\gamma\colon I\to M\mid \gamma(0)=\gamma(1)=x_0\}
\end{align*}
the space of based loops in $M$. We denote by $x_0$ the constant path at $x_0$.

We use the convention that for topological spaces $X,Y\subseteq Z$, we write $C_*(X,Y):=C_*(X,X\cap Y)$ for the relative chains, leaving the intersection implicit.

\addtocontents{toc}{\SkipTocEntry}
\subsection*{Acknowledgements.}
I want to thank my supervisor Nathalie Wahl for suggesting this topic and many enlightening conversation. I also want to thank Isaac Moselle for his helpful comments about the introduction and Section \ref{sec:Sn/G}. I was supported by the Danish National Research Foundation through the Copenhagen Centre for Geometry and Topology (DNRF151). 

\tableofcontents
\section{Various Coproducts}\label{sec:various coproducts}

In this section, we give several definitions of coproducts on the based loop space.
\subsection{The Trivial Coproduct}
We start with the so-called trivial coproduct on $C_*(\Omega M)$.\\

We fix $U$ a contractible neighbourhood of $x_0$ and a representative $\tau\in C^n(U,\{x_0\}^c)$ of the local orientation. We define the map
\begin{align*}
    e\colon \Omega M&\to M,\\
    \gamma&\mapsto\gamma(\tfrac{1}{2}).
\end{align*}
\begin{Def}
    The \textit{Figure 8 space} is defined as
    \begin{align*}
        \mathcal{F}:=\mathcal{F}M:=e^{-1}(x_0)=\{\gamma\in \Omega M\mid \gamma(\tfrac{1}{2})=x_0\}.
    \end{align*}
    The \textit{cutting map} is given by
    \begin{align*}
        c\colon \mathcal{F}&\to \Omega M\times \Omega M,\\
        \gamma&\mapsto (\gamma|_{[0,\frac{1}{2}]},\gamma|_{[\frac{1}{2},1]}).
    \end{align*}
\end{Def}

The inclusion $C_*(e^{-1}(U),\mathcal{F}^c)\to C_*(\Omega M ,\mathcal{F}^c)$ is a homotopy equivalence by excision. By \cite[Proposition 2.21]{hatcher2002algebraic}, we can fix a map
\begin{align*}
    \rho \colon C_*(\Omega M ,\mathcal{F}^c)\to C_*(e^{-1}(U),\mathcal{F}^c)
\end{align*}
that is an inverse on homology.

We fix a homotopy
\begin{align*}
    R_U\colon I\times U\to U
\end{align*}
between $\id_U$ and the constant map $x_0$ that fixes $x_0$.
\begin{Rem}
    In \cite{hingston2010loop}, this is induced by the Riemannian metric on $M$. For a generalization to the topological case, we do not assume that $M$ has a Riemannian metric. This in turn means, we have to fix the above $R_U$.
\end{Rem}
The map
\begin{align*}
    R\colon e^{-1}(U)&\to \mathcal{F},\\
    \gamma&\mapsto \begin{cases}
        \gamma(2t) &\text{if } 0\leq t\leq \frac{1}{4},\\
        R_U(4t-1,\gamma(\frac{1}{2})))& \text{if } \frac{1}{4}\leq t\leq \frac{1}{2},\\
        R_U(3-4t,\gamma(\frac{1}{2})) &\text{if }  \frac{1}{2}\leq t\leq \frac{3}{4},\\
        \gamma(2t-1) &\text{if } \frac{3}{4}\leq t\leq 1
    \end{cases}
\end{align*}
defines a homotopy inverse to the inclusion $\mathcal{F}\to e^{-1}(U)$.

\begin{Def}\label{def:triv coprod}
    The \textit{trivial coproduct} is defined as the composition
    \begin{align*}
        \lor_{\mathrm{triv}}\colon C_*(\Omega M)&\to C_*(\Omega M ,\mathcal{F}^c)\overset{\rho}{\to }C_*(e^{-1}(U),\mathcal{F}^c)\\
        &\overset{e^*\tau\cap}{\to}C_{-n+*}(e^{-1}(U))\overset{R_*}{\to} C_{-n+*}(\mathcal{F})\overset{c_*}{\to}C_{-n+*}(\Omega M\times\Omega M).
    \end{align*}
\end{Def}

\begin{Rem}
    This is the restriction to the based loop space of (a chain level model of) the loop coproduct studied by Tamanoi in \cite{tamanoi2010loop}. Tamanoi shows that the coproduct is zero except on $H_n(M)\subseteq H_*(\Lambda M)$. In particular, the trivial coproduct is zero on homology. However, the coproduct on chains is still of interest as it gives rise to the following more interesting coproducts.
\end{Rem}

\subsection{The Goresky-Hingston Coproduct}\label{subsec:GH coprod}
We next turn our attention to the Goresky-Hingston coproduct as introduced in \cite{goresky2009loop}. We follow the description as in \cite{hingston2017product}.\\

We define the reparametrization map
\begin{align*}
    r\colon I\times \Omega M&\to \Omega M,\\
    (s,\gamma)&\mapsto \left(t\mapsto\begin{cases}
        \gamma(2st) &\text{if } 0\leq t\leq \frac{1}{2};\\
        \gamma(2(1-s)t-1+2s)& \text{if }\frac{1}{2}\leq t\leq 1. 
    \end{cases}\right)
\end{align*}
that parametrizes $\gamma$ such that $r(\gamma)(\frac{1}{2})=\gamma(s)$.

We denote $\mathcal{H}$ as the space of half-constant loops 
\begin{align*}
    \mathcal{H}:=\{\gamma\in \Omega\mid \gamma|_{[0,\frac{1}{2}]}\equiv x_0\text{ or }\gamma|_{[\frac{1}{2},1]}\equiv x_0\}.
\end{align*}
We note that $r$ sends $\partial I\times \Omega M$ and $I\times \{x_0\}$ to $\mathcal{H}$ and $\lor_{\mathrm{triv}}$ sends $C_*(\mathcal{H})$ to $C_{-n+*}(\Omega M\times \{x_0\}\cup \{x_0\}\times \Omega M)$. We therefore get the following definition

\begin{Def}
    The \textit{Goresky-Hingston coproduct} is defined as the composition
    \begin{align*}
        \lor_{\mathrm{GH}}\colon C_*(\Omega M,x_0)&\overset{I\times}{\to}C_{1+*}(I\times \Omega M,\partial I\times \Omega M\cup I\times \{x_0\})\\
        &\overset{r_*}{\to}C_{1+*}(\Omega M,\mathcal{H})\overset{\lor_{\mathrm{triv}}}{\to}C_{1-n+*}((\Omega M,x_0)\times (\Omega M,x_0)).
    \end{align*}
\end{Def}

\begin{Rem}
    This definition is indeed a chain-level model of the original definition as in \cite{goresky2009loop} (see \cite[Theorem 2.13]{hingston2017product}).
\end{Rem}

\subsection{The Avoiding-Stick Definition}\label{subsec:a-s def}
In this subsection, we define a coproduct using sticks avoiding $x_0$. This construction avoids self-intersections coming from $\partial I\times \Omega M$ and thus allows us to lift this coproduct to a coproduct on $C_*(\Omega M)$.\\

We fix $0<\varepsilon<\frac{1}{6}$. 
\begin{Def}
    Let $x_1,\dots,x_k\in M$. A \textit{stick avoiding $x_0,x_1,\dots,x_k$} is a path $\sigma \colon [0,\varepsilon]\to M$ with $\sigma(0)=x_0$ and there exists a chart $\varphi \colon \R^n\cong U\subseteq M$ such that $\{x_0,x_1,\dots, x_k\}\cap U=\{x_0\}$ and $\sigma(t)=\varphi(t,0,\dots,0)$.
\end{Def}

For a stick $\sigma$ avoiding $x_0,x_1,\dots,x_k$, we consider the subspace 
\begin{align}\label{def:Omega sigma}
    \Omega^\sigma M:=\{\gamma\in \Omega \mid \gamma|_{[0,\varepsilon]}=\sigma, \gamma|_{[1-\varepsilon,1]}=\overline{\sigma}\}\subseteq \Omega M
\end{align}
that look like $\sigma$ at the start and at the end. Here $\overline{\sigma}$ denotes the inverse path to $\sigma$.

\begin{Lem}\label{lem:proj to Omega sigma}
    Let $\sigma \colon [0,\varepsilon]\to M$ be a stick avoiding $x_0$. Then the inclusion $\iota \colon \Omega^\sigma M\hookrightarrow\Omega M$ has homotopy inverse
    \begin{align*}
        \pi_\sigma \colon \Omega M & \to \Omega^\sigma M,\\
        \gamma &\mapsto \sigma \cdot\overline{\sigma}\cdot\gamma\cdot  \sigma \cdot\overline{\sigma} 
    \end{align*}
    where the composition is parametrized such that $\gamma$ is followed on $[{\varepsilon},{\varepsilon}]$ in $\frac{1}{1-4\varepsilon}$-speed. 
\end{Lem}
\begin{proof}
    We consider the homotopy
    \begin{align*}
        H(s,\gamma)=\sigma|_{[0,s\varepsilon]} \cdot\overline{\sigma|_{[0,s\varepsilon]}}\cdot\gamma \cdot\sigma|_{[0,s\varepsilon]} \cdot\overline{\sigma|_{[0,s\varepsilon]}} 
    \end{align*}
    where the composition is parametrized such that $\gamma$ is followed on $[2s\varepsilon,1-2s\varepsilon]$ with $\frac{1}{1-4s\varepsilon}$ speed. This defines a homotopy between from $\id_{\Omega M}$ to $\iota\circ \pi_\sigma$.
    
    On the other hand, we can write any element of $\Omega^\sigma M$ uniquely as $\sigma\cdot\beta \cdot\overline\sigma$ for some $\beta\in \Omega_{\sigma(\varepsilon)}M$. Then
    \begin{align*}
        G(s,\sigma\cdot\beta \cdot\overline\sigma)=\sigma \cdot\overline{\sigma|_{[0,s\varepsilon]}}\cdot\sigma|_{[0,s\varepsilon]}\cdot\beta \cdot\overline{\sigma|_{[0,s\varepsilon]}} \cdot\sigma|_{[0,s\varepsilon]} \cdot\overline{\sigma} 
    \end{align*}
    is a homotopy between $\id_{\Omega^\sigma M}$ and $\pi_\sigma\circ \iota$.
\end{proof}

We therefore have a homotopic space $\Omega^\sigma M\simeq \Omega M$. This space has the advantage that we can get rid off the half-constant loops $\mathcal{H}$ because they are isolated from the non-trivial self-intersections in $\mathcal{F}\cap (I\times\Omega^\sigma M)$. Therefore we can apply excision to find that we can only study times $(\varepsilon,1-\varepsilon)$ to detect all non-trivial self-intersections. In precise terms, we have the following key lemma:

\begin{Lem}\label{lem:rho sigma}
    Let $\sigma$ be a stick avoiding $x_0$. The inclusion 
    \begin{align*}
        C_*(r((\varepsilon,1-\varepsilon)\times\Omega^\sigma M),\mathcal{F}^c)\to C_*(r(I\times \Omega^\sigma M),\mathcal{F}^c\cup \mathcal{H})
    \end{align*}
    is a quasi-isomorphism.
\end{Lem}
\begin{proof}
    We want to apply excision. We thus must check that the closure of 
    \begin{align*}
        A:=r(I\times \Omega^\sigma M)\setminus r((\varepsilon,1-\varepsilon)\times\Omega^\sigma M)
    \end{align*}
    is contained in the interior of $\mathcal{F}^c\cup \mathcal{H}$.
    
    It thus suffices to show the following:
    \begin{enumerate}[(1)]
        \item\label{item:A is closed} $A$ is closed;
        \item\label{item:Fc H is open} $\mathcal{F}^c\cup \mathcal{H}$ is open;
        \item\label{item:A in Fc H} $A\subseteq \mathcal{F}^c\cup \mathcal{H}$.
    \end{enumerate}
    For \eqref{item:A is closed}, we check that $r|_{I\times \Omega^\sigma M}$ is a homeomorphism onto its image. We first consider $(s,\gamma)\in I\times \Omega^\sigma M$ with $s\notin \partial I$. We note that $r(s,\gamma)(t)=\sigma(2st)\neq x_0$ for $t\in (0,\tfrac{\varepsilon}{2s})$ and $r(s,\gamma)(t)=\sigma(2(1-s)t-1+2s)\neq x_0$ for $t\in (\frac{1-\varepsilon}{2(1-s)},1)$. This shows that $r(t,\gamma)\notin \mathcal{H}$.
    
    Moreover if $s'\in I$ is such that
    \begin{align*}
        r(s,\gamma)(2s't)=\sigma(t)
    \end{align*}
    for $t\in [0,\varepsilon]$, we have $s=s'$. Therefore we can recover $s$ from $r(s,\gamma)$ and we similarly recover
    \begin{align*}
        \gamma(t)=\begin{cases}
            r(s,\gamma)(\frac{t}{2s}) &\text{if } 0\leq t\leq s;\\
            r(s,\gamma)(\frac{t+1-2s}{2(1-s)}) &\text{if }s\leq t\leq 1.
        \end{cases}
    \end{align*}
    On the other hand if $r(s,\gamma)\in \mathcal{H}$, we find either $r(s,\gamma)|_{[0,\frac{1}{2}]}\equiv x_0$ or $r(s,\gamma)|_{[\frac{1}{2},1]}\equiv x_0$. Both halves are not constant because otherwise $\gamma$ would have been constant but $\gamma$ is not in $\Omega^\sigma M$. If $r(s,\gamma)|_{[0,\frac{1}{2}]}\equiv x_0$, we recover $s=0$ and $\gamma(t)=r(s,\gamma)|_{[\frac{1}{2},1]}(2t-1)$ and, if $r(s,\gamma)|_{[\frac{1}{2},1]}\equiv x_0$, we recover $\gamma(t)=r(s,\gamma)|_{[1,\frac{1}{2}]}(2t)$.
    
    This shows that $r|_{I\times \Omega^\sigma M}$ is a homeomorphism onto its image and thus 
    \begin{align*}
        A=r(I\times \Omega^\sigma M\setminus (\varepsilon,1-\varepsilon)\times\Omega^\sigma M)=r(([0,\varepsilon]\cup [1-\varepsilon,1])\times \Omega^\sigma M).
    \end{align*}
    Therefore $A$ is closed because $([0,\varepsilon]\cup [1-\varepsilon,1])\times \Omega^\sigma M$ is closed in $I\times \Omega^\sigma M$.
    
    We now check \eqref{item:Fc H is open} that $\mathcal{F}^c\cup \mathcal{H}$ is open. On the one hand, $\mathcal{F}^c$ is an open set. On the other hand, the above discussion shows that $\mathcal{H}=r(\partial I\times \Omega^\sigma M)$. This subspace has an open neighbourhood $r(([0,\varepsilon)\cup (1-\varepsilon,1])\times \Omega^\sigma M)$ which is contained in $\mathcal{H}\cup \mathcal{F}^c$ because $r(((0,\varepsilon)\cup (1-\varepsilon,1))\times \Omega^\sigma M)\subseteq \mathcal{F}^c$. This shows that $\mathcal{F}^c\cup \mathcal{H}$ is open.
    
    We finally have to check \eqref{item:A in Fc H} that $A$ is contained in $e(x_0)^c\cup \mathcal{H}$. This follows from the fact that
    \begin{align*}
        A=r(\partial I\times \Omega^\sigma M) \cup r\left(((0,\varepsilon]\times [1-\varepsilon,1))\times \Omega^\sigma M\right)
    \end{align*}
    and that $r(\partial I\times \Omega^\sigma M)\subseteq \mathcal{H}$ and $r\left(((0,\varepsilon]\times [1-\varepsilon,1))\times \Omega^\sigma M\right)\subseteq \mathcal{F}^c$.
    
    Therefore excision implies that 
    \begin{align*}
        C_*(r(I\times \Omega^\sigma M)\setminus A,\mathcal{F}^c\cup \mathcal{H}\setminus A)\to C_*(r(I\times \Omega^\sigma M),\mathcal{F}^c\cup \mathcal{H})
    \end{align*}
    is a homotopy equivalence. The observation
    \begin{align*}
        r(I\times \Omega^\sigma M)\setminus A&=r((\varepsilon,1-\varepsilon)\times\Omega^\sigma M),\\
        \mathcal{F}^c\cup \mathcal{H}\setminus A&\subseteq \mathcal{F}^c
    \end{align*}
    concludes the proof.
\end{proof}

We fix a map
\begin{align*}
    \rho_\sigma\colon C_*(r(I\times \Omega^\sigma M),\mathcal{F}^c\cup \mathcal{H})\to C_*(r((\varepsilon,1-\varepsilon)\times\Omega^\sigma M),\mathcal{F}^c)
\end{align*}
which is an inverse on homology to the above inclusion (e.g.\ see \cite[Proposition 2.21]{hatcher2002algebraic} for a construction of such a map). 

\begin{Def}\label{def:a-s coprod}
    Let $\sigma$ be a stick avoiding $x_0$. The \textit{avoiding-stick coproduct} is defined as the composition
    \begin{align*}
        \lor_\sigma\colon C_*(\Omega M)&\overset{(\pi_\sigma)_*}{\to }C_*(\Omega^\sigma M)\overset{I\times}{\to} C_{1+*}(I\times \Omega^\sigma M, \partial I\times \Omega^\sigma M)\\
        &\overset{r_*}{\to} C_{1+*}(r(I\times \Omega^\sigma M),\mathcal{F}^c\cup \mathcal{H})\\
        &\overset{\rho_\sigma}{\to} C_{1+*}(r((\varepsilon,1-\varepsilon)\times\Omega^\sigma M),\mathcal{F}^c)\\
        &\overset{\lor_{\mathrm{triv}}}{\to}C_{1-n+*}(\Omega M\times \Omega M).
    \end{align*}
\end{Def}

\begin{Rem}
    We note that because the first map in the definition of $\lor_{\mathrm{triv}}$ in Definition \ref{def:triv coprod} is 
    \begin{align*}
        C_*(\Omega M)\to C_*(\Omega M,\mathcal{F}^c)
    \end{align*}
    the avoiding-stick coproduct lands indeed in non-relative chains.
\end{Rem}

The avoiding-stick coproduct is a lift of the Goresky-Hingston coproduct:

\begin{Prop}\label{prop:a-s-coprod lifts GH}
    The map $i\colon (\Omega M,\varnothing)\to (\Omega M,x_0)$ intertwines the Goresky-Hingston coproduct and the avoiding-stick coproduct on homology, that is, the following diagram commutes:
    \begin{center}
        \begin{tikzcd}
            H_*(\Omega M)\ar[d, "\lor_{\sigma}"] \ar[r,"i_*"] & H_*(\Omega M,x_0)\ar[d, "\lor_{\mathrm{GH}}"]\\
            H_*(\Omega M\times \Omega M)\ar[r, "(i\times i)_*"]& H_*((\Omega M,x_0)\times (\Omega M,x_0)).
        \end{tikzcd}
    \end{center}
\end{Prop}
\begin{proof}
    We denote $\iota_\sigma$ for the inclusion $\Omega^\sigma M\to \Omega M$ and its induced maps. We consider the diagram
    \begin{center}
        \begin{tikzcd}
            C_*(\Omega M)\ar[d, "(\pi_\sigma)_*"]  \ar[rd, "i_*"]\\
            C_*(\Omega^\sigma M) \ar[r, "(\iota_\sigma)_*"] \ar[d, "I\times "] & C_*(\Omega M,x_0)\ar[d, "I\times "]\\
            C_{1+*}(I\times \Omega^\sigma M,\partial I \times \Omega^\sigma M) \ar[r, "(\iota_\sigma)_*"] \ar[d, "r_*"] & C_{1+*}(I\times \Omega M,\partial I \times \Omega M\cup I\times \{x_0\})\ar[d, "r_*"]\\
            C_{1+*}(r(I\times \Omega^\sigma M), \mathcal{F}^c\cup \mathcal{H}) \ar[r,"(\iota_\sigma)_*"] \ar[d, "\rho_\sigma"] & C_{1+*}(\Omega M, \mathcal{F}^c\cup \mathcal{H})\ar[dd, "\lor_{\mathrm{triv}}"]\\
            C_{1+*}(r((\varepsilon,1-\varepsilon)\times \Omega^\sigma M),\mathcal{F}^c)\ar[d, "\lor_{\mathrm{triv}}"]\ar[ru, "(\iota_\sigma)_*"]\\
            C_{1-n+*}(\Omega M\times \Omega M)\ar[r, "(i\times i)_*"] \ar[r] & C_{1-n+*}((\Omega M,x_0)\times (\Omega M,x_0)).
        \end{tikzcd}
    \end{center}
    The top triangle commutes in homology because $(\pi_\sigma)_*$ and $(\iota_\sigma)_*$ are homotopy inverses (see Lemma \ref{lem:proj to Omega sigma}). The bottom triangle commutes in homology because $\rho_\sigma$ is a homotopy inverse to an inclusion of chain complexes and inclusions of chain complexes commute. The rest of the diagram commutes because the maps from left to right are come from inclusions of spaces. Therefore the diagram commutes in homology.
    
    Composition along the bottom left computes $(i\times i)_*\circ \lor_\sigma$ and composition along the top left computes $\lor_{\mathrm{GH}}\circ i_*$.
\end{proof}

A natural question to ask is now the following:
\begin{Ques}\label{ques:depend on sigma}
    Does the absolute coproduct depend on the choice of stick $\sigma$?
\end{Ques}

The answer is yes, if $n=1$ (see Example \ref{exmp:S1 part 1}), but no, if $n>1$ (see Corollary \ref{cor:independence of sigma}).

\begin{Exmp}\label{exmp:S1 part 1}
    For $M=S^1$, the lift is not independent of choice of $\sigma$. We consider $S^1\subseteq \C$. We have
    \begin{align*}
        H_*(\Omega S^1)=\Z[x^{\pm1}]
    \end{align*}
    as a Pontryagin algebra with $x$ in degree $0$. Because every group is free, the Künneth map $H_*(\Omega S^1\times \Omega S^1)\cong H_*(\Omega S^1)\otimes H_*(\Omega S^1)$ is an isomorphism. We can therefore consider the avoiding-stick coproduct as a map $H_*(\Omega S^1)\to H_*(\Omega S^1)\otimes H_*(\Omega S^1)$.
    
    We make the choice that our basepoint is $x_0:=1\in \C$. Moreover, we choose 
    \begin{align*}
        \gamma_1\colon I/\partial I&\to S^1,\\
        t&\mapsto \exp(2\pi i t)
    \end{align*}
    as a representative of $x\in H_0(\Omega S^1)$.
    
    We fix $\varepsilon:=\frac{1}{8}$. Up to homotopy, there are two choices for $\sigma$:
    \begin{align*}
        \sigma^+\colon [0,\tfrac{1}{8}]&\to S^1,\\
        t&\mapsto \exp(2\pi i t);\\
        \sigma^-\colon [0,\tfrac{1}{8}]&\to S^1,\\
        t&\mapsto \exp(-2\pi i t).
    \end{align*}
    We set $x^+:=\sigma^+(\tfrac{1}{8})$ and $x^-:=\sigma^-(\tfrac{1}{8})$. We compute the coproduct $\lor^+:=\lor_{\sigma^+}$: 
    \begin{itemize}
        \item We first compute $\lor^+(1)$. The element $1=x^0\in H_*(\Omega S^1)$ is represented by the constant loop. It gets sent to $\sigma^+\cdot\overline{\sigma^+}\cdot\sigma^+\cdot\overline{\sigma^+}\in \Omega^{\sigma^+} S^1$ under the map $\pi_{\sigma^+}$. This loop is homotopic to $\gamma:=\sigma^+\cdot c_{x^+}\cdot\overline{\sigma^+}$ as an element in $\Omega^{\sigma^+}S^1$. Here $c_{x^+}$ denotes the constant loop at $x^+$. This loop has no $t\in [\frac{1}{8},\frac{7}{8}]$ such that $\gamma(t)=0$. Therefore, $r_*(I\times \gamma)\in C_*(\mathcal{F}^c\cup \mathcal{H})$ and we have $\lor^+(1)=0$.

        \item Next, we compute $\lor^+(t^k)$ for $k>0$. This element is represented by 
        \begin{align*}
            \gamma_k(t)= \exp(2\pi i kt).
        \end{align*}
        It gets sent to $\sigma^+\cdot\overline{\sigma^+}\cdot\gamma_k\cdot\sigma^+\cdot\overline{\sigma^+}\in \Omega^{\sigma^+} S^1$ under the map $\pi_{\sigma^+}$. This is homotopic to the map 
        \begin{align*}
            \gamma_k'(t)=\begin{cases}
                \sigma^+(t) &\text{if } 0\leq t\leq \frac{1}{8};\\
                \sigma^+(1-t) &\text{if }\frac{7}{8}\leq t\leq 1;\\
                \exp\left(\frac{2\pi i (t-\frac{1}8{)}}{1-\frac{1}{4}}+\frac{2\pi i}{8}\right) & \text{if }\frac{1}{8}\leq t\leq \frac{7}{8};
            \end{cases}
        \end{align*}
        as an element in $\Omega^{\sigma^+} S^1$.
        
        There are $k$ times $t\in (\frac{1}{8},\frac{7}{8})$ such that $\gamma_k(t)=x_0$: concretely, define $t_l$ such that 
        \begin{align*}
            \frac{2\pi i (t_l-\frac{1}8{)}}{1-\frac{1}{4}}+\frac{2\pi i}{8}=2\pi i l
        \end{align*}
        for $l=1,\dots, k$. These are transversal intersections of the basepoint. Therefore, they get detected by capping with $e^*\tau\in C^0(e^{-1}(U),\mathcal{F}^c)$.
        
        Cutting at $t_k$, we get the loop $\gamma_k$ in the first variable and $\sigma^+\cdot\overline{\sigma^+}\simeq c_{x_0}$ in the second variable. After cutting at $t_{l}$ for $l<k$, we get $\sigma^+\cdot\overline{\sigma^+}\cdot\gamma_l\simeq \gamma_l$ in the first variable and $\gamma_{k-l}\cdot\sigma^+\cdot\overline{\sigma^+}\simeq \gamma_{k-l}$.
        
        We therefore find
        \begin{align*}
            \lor^+(x^k)=\sum_{l=1}^k x^l\otimes x^{k-l}.
        \end{align*}

        \item An analogous computation for $k<0$ shows
        \begin{align*}
            \lor^+(x^k)=-\sum_{l=k+1}^{0}x^l\otimes x^{k-l}.
        \end{align*}
        The negative sign comes from the fact that all intersections pass $x_0$ in negative direction.
    \end{itemize}
    We summarize our results for $\lor^+$ as
    \begin{align*}
        \lor^+(t^k)=\begin{cases}
            0 &\text{if }k=0;\\
            \sum_{l=1}^k x^l\otimes x^{k-l} &\text{if } k>0;\\
            -\sum_{l=k+1}^0 x^l\otimes x^{k-l} &\text{if }k<0.
        \end{cases}
    \end{align*}
    An analogous computation for $\sigma^-$ shows
    \begin{align*}
        \lor^-(x^k)=\begin{cases}
            0 &\text{if }k=0;\\
            \sum_{l=0}^{k-1} x^l\otimes x^{k-l} &\text{if } k>0;\\
            -\sum_{l=k}^{-1} x^l\otimes x^{k-l} &\text{if }k<0.
        \end{cases}
    \end{align*}
\end{Exmp}

\begin{Rem}
    In \cite[Subsection 8.4]{cieliebak2023loop}, Cieliebak, Hingston and Oancea study four different lifts of the Goresky-Hingston coproduct. Two of them correspond to an avoiding-stick coproduct: it holds $\lor^+=\lambda_{v_-,v_-}$ and $\lor^-=\lambda_{v_+,v_+}$ for $\lambda_{v_-,v_-}$ and $\lambda_{v_+,v_+}$ as in \cite[Remark 8.8]{cieliebak2023loop}.

    These coproducts are neither cocommutative or coassociative. In \cite[Proposition 8.6]{cieliebak2023loop}, Cieliebak, Hingston and Oancea also compute two coproducts which on the based loop space give
    \begin{align*}
        \lambda_+(x^k)=\begin{cases}
            \sum_{l=0}^k x^l\otimes x^{k-l} &\text{if }k\geq 0;\\
            -\sum_{l=k+1}^{-1}x^l\otimes x^{k-1} &\text{if }k<0;
        \end{cases}
    \end{align*}
    and
    \begin{align*}
        \lambda_+(x^k)=\begin{cases}
            \sum_{l=1}^{k-1} x^l\otimes x^{k-l} &\text{if }k>0;\\
            -\sum_{l=k}^{0}x^l\otimes x^{k-1} &\text{if }k\leq 0.
        \end{cases}
    \end{align*}
    Such results can be obtained if one replaces $\Omega^\sigma M$ by
    \begin{align*}
        \{\gamma \in \Omega M\mid \gamma|_{[0,\frac{1}{8}]}=\sigma^+, \gamma|_{[\frac{7}{8},1]}=\overline{\sigma^{-}}\}
    \end{align*}
    which gives $\lambda_+$ or by
    \begin{align*}
        \{\gamma \in \Omega M\mid \gamma|_{[0,\frac{1}{8}]}=\sigma^-, \gamma|_{[\frac{7}{8},1]}=\overline{\sigma^{+}}\}
    \end{align*}
    which gives $\lambda_-$.
    
    In general there is no need to restrict to a subspace where $\gamma(t)=\gamma(1-t)=\sigma(t)$ for small $t$ as we do in our definition of the avoiding-path coproduct. One can choose a path for the exit at $0\leq t\leq \varepsilon$ and another path for the entrance at $1-\varepsilon\leq t\leq 1$ and run the same machinery. However if $n>1$, we show in Corollary \ref{cor:independence of sigma} that the coproduct does not depend on the choice of the paths. For convenience, we thus fix the same path for exit and entrance.
\end{Rem}

\section{Cutting at Other Points}\label{sec:Cutting at Other Points}
In this section, we study maps defined by cutting at points other than the basepoint. We relate them to the coproducts we have defined so far. Cutting at other points naturally shows up when we consider maps between manifolds.

\subsection{The Trivial Pre-Coproduct}
We fix some point $x\in M$. We denote 
\begin{align*}
    \mathcal{P}_{x_0\to x}&:=\{\gamma\colon I\to M\mid \gamma(0)=x_0, \gamma(1)=x\},\\
    \mathcal{P}_{x\to x_0}&:=\{\gamma\colon I\to M\mid \gamma(0)=x, \gamma(1)=x_0\},\\
    \mathcal{F}_{x}&:=e^{-1}(x):=\{\gamma \in \Omega M \mid  \gamma(\tfrac{1}{2})=x\}.
\end{align*}
We fix some contractible neighbourhood $U_{x}$ of $x$. We can define
a homotopy inverse 
\begin{align*}
    \rho_x\colon C_*(\Omega M ,\mathcal{F}_x^c)\to C_*(e^{-1}(U_x),\mathcal{F}^c)
\end{align*}
to the inclusion and a map
\begin{align*}
    R_{x}\colon e^{-1}(U_{x})\to \mathcal{F}_{x}
\end{align*}
akin to $R$. Similarly, we can define a cutting map 
\begin{align*}
    c_{x}\colon\mathcal{F}_{x}\to \mathcal{P}_{x_0\to x}\times \mathcal{P}_{x\to x_0}.
\end{align*}
The orientation of $M$ gives us a local orientation $\tau_{x}\in C^n(U_{x},\{x\}^c)$. Using this, we can define an analogue of the trivial coproduct

\begin{Def}
    Let $x\in M$. The \textit{trivial pre-coproduct at $x$} is given as the composition
    \begin{align*}
        \lor_{\mathrm{triv}}^x\colon C_*(\Omega M)&\to C_*(\Omega M ,\mathcal{F}_x^c)\overset{\rho_x}{\to }C_*(e^{-1}(U_x),\mathcal{F}_x^c)\\
        &\overset{e^*\tau_x\cap}{\to}C_{-n+*}(e^{-1}(U_x))\overset{(R_x)_*}{\to} C_{-n+*}(\mathcal{F}_x)\overset{(c_x)_*}{\to}C_{-n+*}(\mathcal{P}_{x_0\to x}\times\mathcal{P}_{x\to x_0}).
    \end{align*}
\end{Def}

\begin{Rem}
    The definition of the trivial pre-coproduct at $x_0$ gives the definition of the trivial coproduct as in Definition \ref{def:triv coprod}.
\end{Rem}

\begin{figure}[H]
    \centering
    \includegraphics[width=0.5\linewidth]{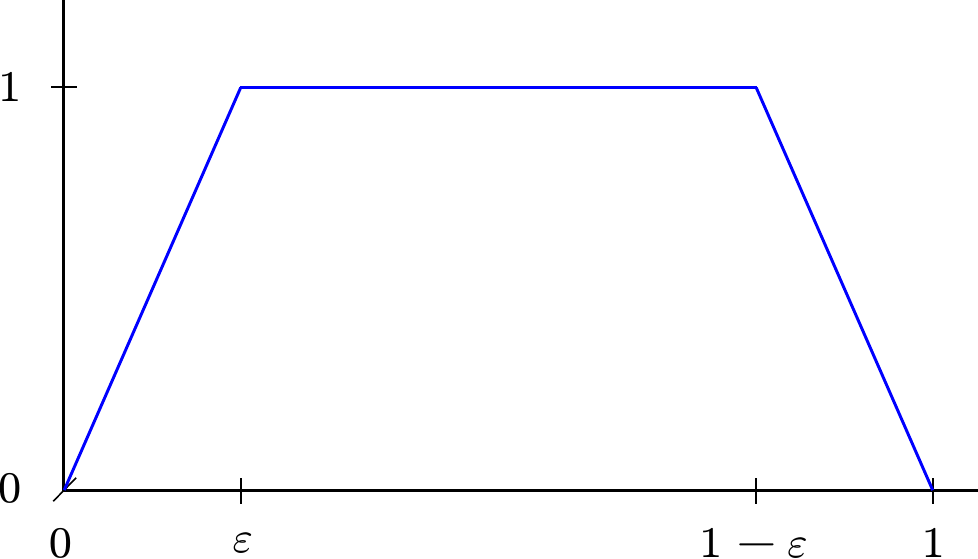}
    \caption{The graph of the function $\chi$.}
    \label{fig:graph of chi}
\end{figure}

We now want to relate different pre-coproducts to each other. First, we note that a path $\tau$ from $x_1$ to $x_2$ defines a homotopy equivalence
\begin{align*}
    r_\tau\times l_{\overline\tau}\colon \mathcal{P}_{x_0\to x_1}\times \mathcal{P}_{x_1\to x_0}&\to \mathcal{P}_{x_0\to x_2}\times \mathcal{P}_{x_2\to x_0},\\
    (\gamma_1,\gamma_2)&\mapsto (\gamma_1\cdot \tau,\overline{\tau}\cdot\gamma_2).
\end{align*}
We denote the simplest piece-wise linear map $\chi\colon I\to I$ such that
\begin{align*}
    \chi(0)&=0,& \chi|_{[\varepsilon,1-\varepsilon]}&\equiv 1, &\chi(1)=0.
\end{align*}
We note that $\chi(t)=\chi(1-t)$. (See Figure \ref{fig:graph of chi}.)

\begin{Lem}\label{lem:triv commuting with homotopy}
    Let $x_1,x_2\in M$. Let $\phi^s\colon M\to M$ for $s\in [0,1]$ be a $1$-parameter family of homeomorphisms such that $\phi^1(x_1)=x_2$.
    
    Let $\tau$ be the path $\tau(t)=\phi^t(x_1)$ and 
    \begin{align*}
        \Phi\colon \Omega M&\to  \Omega M,\\
        \gamma&\mapsto (t\mapsto \phi^{\chi(t)}\circ \gamma(t)).
    \end{align*}
    It holds
    \begin{align*}
        (r_\tau\times l_{\overline{\tau}})_*\circ\lor_{\mathrm{triv}}^{x_1}=\lor_{\mathrm{triv}}^{x_2}\circ \Phi_*
    \end{align*}
    as maps $H_*(\Omega M,\mathcal{F}_{x_1}^c)\to H_{-n+*}(\mathcal{P}_{x_0\to x_2}\times\mathcal{P}_{x_2\to x_0})$.
\end{Lem}
\begin{proof}
    The map $\Phi_*$ sends $C_*(\Omega M,\mathcal{F}_{x_1}^c)$ to $C_*(\Omega M,\mathcal{F}_{x_2}^c)$. Indeed, $\Phi$ sends $\mathcal{F}_{x_1}^c$ to $\mathcal{F}_{x_2}^c$. Moreover in homology, the maps $\rho_{x_i}$ and $R_{x_i}$ can be defined such that they are intertwined by $\Phi$.
    
    The naturality of the cap product gives that the diagram
    \begin{center}
        \begin{tikzcd}
            C_*(e^{-1}(U_{x_1}),e(x_1)^c)\ar[r, "\Phi_*"]\ar[d, "\Phi^*e^*\tau_{x_2}\cap"] & C_*(e^{-1}(U_{x_2}),e(x_2)^c) \ar[d, "e^*\tau\cap"]\\
            C_{-n+*}(e^{-1}(U_{x_1}))\ar[r, "\Phi_*"] &C_{-n+*}(e^{-1}(U_{x_2})).
        \end{tikzcd}
    \end{center}
    We set $\phi:=\phi^1$ and find in homology
    \begin{align*}
        \Phi^*e^*\tau_{x_2}=e^*\phi^*\tau_{x_2}=\deg(\phi)_{x_1}e^*\tau_{x_1}.
    \end{align*}
    Because $\phi$ is homotopic to the identity, we have that $\deg(\phi)_{x_1}=1$. Therefore $\Phi$ intertwines the map $e^*\tau_{x_i}\cap$.
    
    We have shown so far that $\Phi_*$ intertwines all but the last map in the definition of the trivial pre-coproduct. We thus have:
    \begin{align*}
        \lor_{\mathrm{triv}}^{x_2}\circ \Phi_*=(c_{x_2})_*\circ \Phi_*\circ (R_{x_1})_*\circ (e^*\tau_{x_1}\cap)\circ \rho_{x_1}
    \end{align*}
    as maps $H_*(\Omega M,\mathcal{F}_{x_1}^c)\to H_{-n+*}(\mathcal{P}_{x_0\to x_2}\times\mathcal{P}_{x_2\to x_0})$.
    
    Finally we note that $\Phi\colon \mathcal{F}_{x_1}\to \mathcal{F}_{x_2}$ is homotopic to the map
    \begin{align*}
        \mathcal{F}_{x_1}&\to \mathcal{F}_{x_2},\\
        \gamma&\mapsto \gamma|_{[0,\frac{1}{2}]}\cdot \tau\cdot \overline{\tau}\cdot \gamma
    \end{align*}
    through the homotopies 
    \begin{align*}
        [\gamma\mapsto\Phi(\gamma)]&\simeq \left[\gamma\mapsto \Phi(\gamma|_{[0,\frac{1}{2}]}\cdot c_{x_1}\cdot c_{x_1}\cdot \gamma|_{[\frac{1}{2},1]} )=\left(t\mapsto\begin{cases}
            \phi^{\chi(2t)}(\gamma(2t)) &\text{if }0\leq t\leq \frac{1}{4};\\
            \phi^{1}(x_1) &\text{if }\frac{1}{4}\leq t\leq \frac{3}{4};\\
            \phi^{\chi(2t-1)}(\gamma(2t-1)) &\text{if }\frac{3}{4}\leq t\leq 1;
        \end{cases}\right)\right]\\
        &\simeq\left[\gamma \mapsto\left(t\mapsto\begin{cases}
            \phi^0(\gamma(2t)) &\text{if } 0\leq t\leq \frac{1}{4};\\
            \phi^{4t-1}(x_1) &\text{if }\frac{1}{4}\leq t\leq\frac{1}{2};\\
            \phi^{3-4t}(x_1) &\text{if }\frac{1}{2}\leq t\leq \frac{3}{4};\\
            \phi^0(\gamma(2t-1)) &\text{if }\frac{3}{4}\leq t\leq 1;
        \end{cases}\right)\right]\\
        &=[\gamma\mapsto \gamma|_{[0,\frac{1}{2}]}\cdot \tau\cdot \overline{\tau}\cdot \gamma].
    \end{align*}
    Therefore, $c_{x_2}\circ \Phi$ is homotopic to $(r_\tau\times l_{\overline{\tau}})\circ c_{x_1}$ and, on homology, we have
    \begin{align*}
        (r_\tau\times l_{\overline{\tau}})\circ \lor_{\mathrm{triv}}^{x_1}=\lor_{\sigma_2,x_2}\circ \Phi_*.
    \end{align*}
\end{proof}

\subsection{The Avoiding-Stick Pre-Coproduct}
In this subsection, we define and study an analogous construction to the avoiding-stick coproduct where we cut at different points. We relate the maps that arise from cutting at different points to each other. Using this, we can prove a $\pi_1$-invariance result in Corollary \ref{cor:pi1 invariance} and prove the Sullivan relation in Proposition \ref{prop:sullivan} between the avoiding-stick coproduct and the Pontryagin product.\\

Analogous to Lemma \ref{lem:rho sigma}, one proves the following lemma:

\begin{Lem}
    Let $\sigma$ be a stick avoiding $x_0,x$. The inclusion 
    \begin{align*}
        C_*(r((\varepsilon,1-\varepsilon)\times\Omega^\sigma M),\mathcal{F}_x^c)\to C_*(r(I\times \Omega^\sigma M),\mathcal{F}_x^c\cup \mathcal{H})
    \end{align*}
    is a homotopy equivalence by excision.
\end{Lem}

We can thus fix a homotopy inverse
\begin{align*}
    \rho^x_\sigma\colon C_*(r(I\times \Omega^\sigma M),\mathcal{F}_x^c\cup \mathcal{H})\to C_*(r((\varepsilon,1-\varepsilon)\times\Omega^\sigma M),\mathcal{F}_x^c).
\end{align*}

\begin{Def}\label{def:a-s pre-coprod}
    Let $\sigma$ be a stick avoiding $x_0,x$. The \textit{avoiding-stick pre-coproduct at $x$} is defined as the composition
    \begin{align*}
        \lor_\sigma^x\colon C_*(\Omega M)&\overset{(\pi_\sigma)_*}{\to }C_*(\Omega^\sigma M)\overset{I\times}{\to} C_{1+*}(I\times \Omega^\sigma M, \partial I\times \Omega^\sigma M)\\
        &\overset{r_*}{\to} C_{1+*}(r(I\times \Omega^\sigma M),\mathcal{F}_x^c\cup \mathcal{H})\\
        &\overset{\rho_\sigma^x}{\to} C_{1+*}(r((\varepsilon,1-\varepsilon)\times\Omega^\sigma M),\mathcal{F}_x^c)\\
        &\overset{\lor_{\mathrm{triv}}^x}{\to}C_{1-n+*}(\mathcal{P}_{x_0\to x}\times\mathcal{P}_{x\to x_0}).
    \end{align*}
\end{Def}

\begin{Rem}
    The avoiding-stick pre-coproduct at $x_0$ is the same as the avoiding-stick coproduct as in Definition \ref{def:a-s coprod}.
\end{Rem}

\begin{figure}[H]
    \centering
    \includegraphics[width=0.35\linewidth]{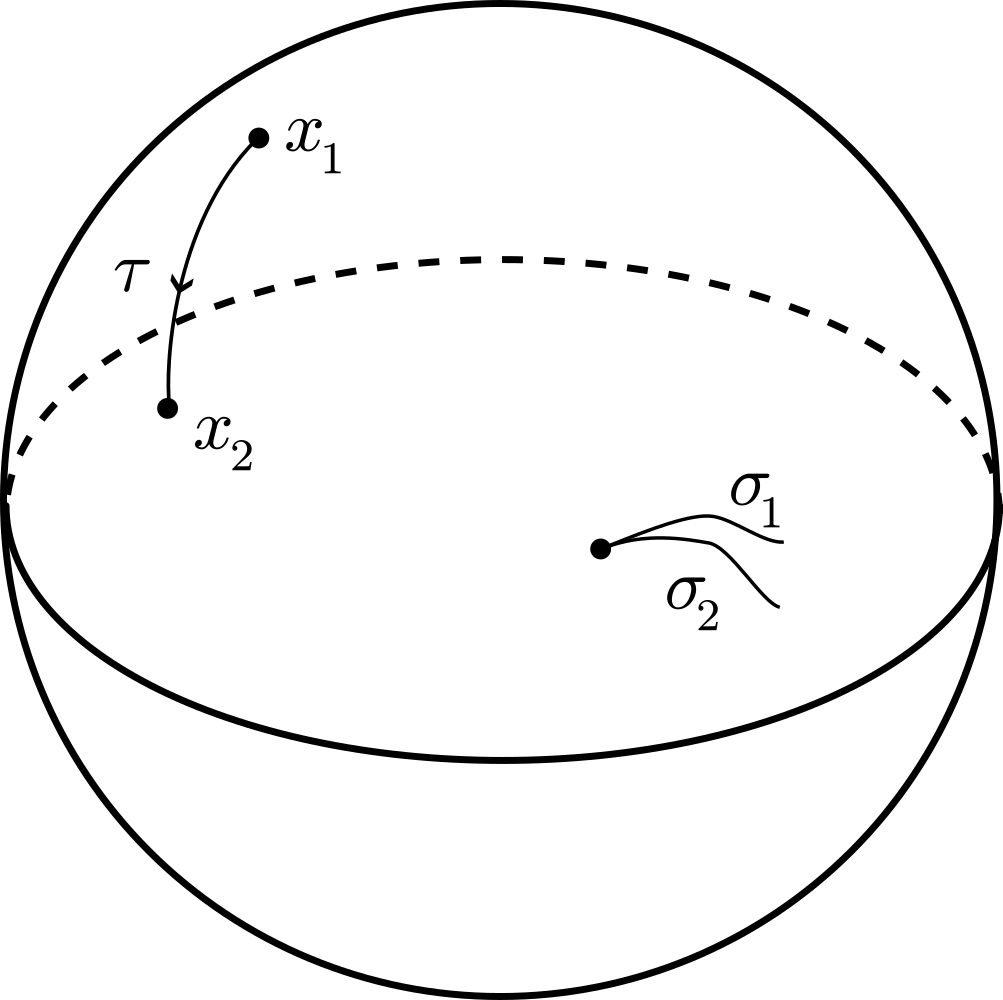}
    \caption{The set up of Proposition \ref{prop:a-s commuting with homotopy} for $M=S^2$ where the $1$-parameter family is given by a flow from the north pole to the south pole. The path $\sigma_2$ is then given by a \enquote{bent version} of $\sigma_2$.}
    \label{fig:1-parameter family}
\end{figure}

\begin{Prop}\label{prop:a-s commuting with homotopy}
    Let $x_1,x_2\in M$ and $\sigma_1$ be a stick avoiding $x_0,x_1$. Let $\phi^s\colon M\to M$ be a $1$-parameter family of homeomorphisms (see Figure \ref{fig:1-parameter family}) such that
    \begin{enumerate}[(i)]
        \item $\phi^0=\id$,
        \item $\phi^1(x_1)=x_2$,
        \item $\sigma_2(t):=\phi^{\chi(t)}\circ \sigma_1(t)\neq x_2$ is a stick avoiding $x_0,x_2$.
    \end{enumerate}
    Let $\tau$ be the path $\tau(t)=\phi^t(x_1)$, then in homology the following formula holds
    \begin{align*}
        (r_\tau\times l_{\overline{\tau}})\circ\lor_{\sigma_1}^{x_1}=\lor_{\sigma_2}^{x_2}
    \end{align*}
    as maps from $H_*(\Omega M)\to H_*(\mathcal{P}_{x_0\to x_2}\times\mathcal{P}_{x_2\to x_0})$.
\end{Prop}
\begin{proof}
    We again denote the map
    \begin{align*}
        \Phi\colon \Omega M&\to  \Omega M,\\
        \gamma&\mapsto (t\mapsto \phi^{\chi(t)}\circ \gamma(t)).
    \end{align*}
    By Lemma \ref{lem:triv commuting with homotopy}, we have $\lor^{x_2}\circ \Phi_*=(r_\tau\times l_{\overline{\tau}})_*\circ\lor_{\mathrm{triv}}^{x_1}$. We therefore need to show that $\Phi_*$ intertwines the other maps in Definition \ref{def:a-s pre-coprod}.
    
    We note that $\id\times \Phi$ sends $(I\times \Omega^{\sigma_1} M,\partial I\times \Omega^{\sigma_1}M)$ to $(I\times \Omega^{\sigma_2} M,\partial I\times \Omega^{\sigma_2}M)$. In fact, we have $\Phi\circ \pi_{\sigma_1}=\pi_{\sigma_2}\circ \Phi$.
    
    It thus remains to show that on homology it holds:
    \begin{align}\label{eq:Phi commutes with r and rho}
        \Phi_*\circ \rho_{\sigma_1}\circ r_*=\rho_{\sigma_2}\circ r_*\circ (\id\times \Phi)_*
    \end{align}
    as maps $H_*(I\times \Omega^{\sigma_1} M,\partial I\times \Omega^{\sigma_1})\to H_*(\Omega^{\sigma_2}M, e^{-1}(x_2))$.
    
    In the proof of Lemma \ref{lem:rho sigma}, we showed that $r$ is a homeomorphism onto its image. We define ${}^r\Phi:=r\circ (\id\times \Phi)\circ r^{-1}$. By construction, we thus have $r\circ (\id\times\Phi)={}^r\Phi\circ r$. It remains to show that $\Phi_*\circ \rho_{x_1}=\rho_{x_2}\circ {}^r\Phi_*$ on homology. The maps $\rho_{\sigma_i}$ were constructed as inverses to the inclusions
    \begin{align*}
        C_{*}(r((\varepsilon,1-\varepsilon)\times\Omega^{\sigma_i} M),\mathcal{F}_{x_i}^c)\to C_{*}(r(I\times \Omega^{\sigma_i} M),\mathcal{F}_{x_i}^c\cup \mathcal{H})
    \end{align*}
    These inclusions are intertwined by ${}^r\Phi_*$ and thus so are $\rho_{\sigma_i}$ on homology.
    
    Finally, the square
    \begin{center}
        \begin{tikzcd}
            C_{*}(r((\varepsilon,1-\varepsilon)\times\Omega^{\sigma_1} M),\mathcal{F}_{x_1}^c) \ar[r, "({}^r\Phi)_*"] \ar[d] & C_{*}(r((\varepsilon,1-\varepsilon)\times\Omega^{\sigma_2} M),\mathcal{F}_{x_2}^c)\ar[d]\\
            C_{*}(\Omega M,\mathcal{F}_{x_1}^c) \ar[r, "\Phi_*"] &C_{*}(\Omega M,\mathcal{F}_{x_2}^c)
        \end{tikzcd}
    \end{center}
    commutes in homology because ${}^r\Phi(\gamma)$ and $\Phi(\gamma)$ are the same loop up to reparametrization.
    
    This shows \eqref{eq:Phi commutes with r and rho} and thus concludes the proof.
\end{proof}

This result proves to be very powerful as we can often find a $1$-parameter family $\phi^s$ satisfying the assumptions and achieving some other conditions.

\begin{Cor}\label{cor:independence of sigma}
    Let $\sigma_1,\sigma_2$ be two paths avoiding $x_0,x_1$ in a manifold of dimension $n>1$. The two pre-coproducts $\lor_{\sigma_1,x_1}$ and $\lor_{\sigma_2,x_1}$ are the same in homology.
    
    In particular for $x_0=x_1$, we have that $\lor_{\sigma_1}=\lor_{\sigma_2}$ in homology, that is that the avoiding-stick coproduct is independent of the choice of the avoiding-stick in homology.
\end{Cor}
\begin{proof}
    Our strategy to prove this corollary is to find a $1$-parameter family $\phi^t$ that satisfies the assumptions for Proposition \ref{prop:a-s commuting with homotopy} such that $\tau(t)=\phi^t(x_1)$ is a contractible loop based at $x_1$.
    
    By assumption, we can write, for $i=1,2$, $\sigma_i(t)=\varphi_i(t,0,\dots,0)$ for charts $\varphi_i\colon \R^n\cong U_i\subseteq M$.
    
    If $x_1\neq x_0$, we have that $x_1\notin U_1,U_2$. Moreover, by applying Proposition \ref{prop:a-s commuting with homotopy} to $\sigma_2$ and a $1$-parameter family that shrinks $U_2$, we can assume that $U_2\subseteq U_1$.
    
    We consider $\sigma_1(s\varepsilon)$ as an isotopy of topological embeddings (historically called locally flat maps) of the point $x_0$ into $M$. By the isotopy extension theorem for topological manifolds \cite[Corollary 1.4]{edwards1971deformations}, we can find an isotopy $H^s_1\colon M\to M$ such that $H_1^s(x_0)=\sigma_1(s\varepsilon)$ that is supported in $U_1$.
    
    Because $U_1$ is contractible, we can find a path $\rho$ between $y:=H^{-1}_1(x_0)$ and $x_0$. We then can extend the isotopy
    \begin{align*}
        h^s_2\colon \{x_0,y\}&\to M,\\
        x_0&\mapsto \sigma_2(\varepsilon s),\\
        y&\mapsto \rho(s)
    \end{align*}
    to a isotopy $H_2^s\colon M\to M$.
    
    The $1$-parameter family $\phi^s:=H^{s}_2\circ H_1^{-s}$ satisfies
    \begin{enumerate}[(i)]
        \item $\phi^0=H^0_2\circ H^0_1=\id$;
        \item $\phi^1(x_1)=H_2^{1}\circ H_1^{-1}(x_1)=x_1$, if $x_1\neq x_0$, or 
        \begin{align*}
            \phi^1(x_1)=\phi^1(x_0)=H^1_2(H^{-1}_1(x_0))=H_2^1(y)=\rho(1)=x_0=x_1,
        \end{align*}
        if $x_1=x_0$;
        \item for $0\leq t\leq \varepsilon$, we set $s:=t/\varepsilon\in I$ and compute
        \begin{align*}
            \phi^{\chi(t)}(\sigma_1(t))=\phi^{s}(\sigma_1(s\varepsilon))=H_2^s(H_1^{-s}(\sigma_1(s\varepsilon))=H_2^s(x_0)=\sigma_2(s\varepsilon)=\sigma_2(t).
        \end{align*}
    \end{enumerate}
    We therefore can apply Proposition \ref{prop:a-s commuting with homotopy} and find in homology that $(r_\tau\times l_{\overline{\tau}})\circ \lor_{\sigma_1}^{x_1}=\lor_{\sigma_2}^{x_1}$ for $\tau(t)=\phi^t(x_1)$.
    
    If $x_1\neq x_0$, we have $\phi^t(x_1)=H_2^t(H_1^{-t}(x_1))\equiv x_1$. If $x_1=x_0$, then $\phi^t(x_1)$ is a loop in $U_1$ and thus contractible. In both cases, we deduce $(r_\tau\times l_{\overline{\tau}})=\id$ in homology which shows the Corollary.
\end{proof}

This in particular answers Question \ref{ques:depend on sigma} whether or not the avoiding-stick coproduct depends on the stick $\sigma$. We are therefore justified to use the following notation:

\begin{Not}
    For $n>1$, we denote the avoiding-stick coproduct
    \begin{align*}
        \lor_\sigma\colon  H_*(\Omega M)\to H_{1-n+*}(\Omega M\times \Omega M)
    \end{align*}
    by $\lor$ or $\lor_{\Omega M}$. Corollary \ref{cor:independence of sigma} for $x_0=x_1$ shows that this is well-defined.
    
    Similarly, we denote for $x\in M$ and $\sigma$ a stick avoiding $x_0,x$, $\lor^x:=\lor_{\sigma}^x$.
\end{Not}

\begin{Cor}\label{cor:relate pre-coprod for n>1}
    Assume that the dimension is $n>1$. For any path $\tau$ from $x_1$ to $x_2$, it holds in homology that 
    \begin{align*}
        (r_\tau\times l_{\overline{\tau}})\circ \lor^{x_1}=\lor^{x_2}.
    \end{align*}
\end{Cor}
\begin{proof}
    We again want to apply Proposition \ref{prop:a-s commuting with homotopy}.
    
    We note that for $\tau_1$ a path from $x_1$ and $x_2$ and $\tau_2$ a path from $x_2$ to $x_3$, the statement for $\tau_1$ and $\tau_2$ implies it for $\tau_1\cdot \tau_2$. Therefore, we can decompose $\tau$ and thus assume that $\tau$ is contained in a contractible neighbourhood. Moreover, we note that the statement only depends on the homotopy type (relative to $\partial I$) of $\tau$. We can thus assume that $\tau$ is the straight line from $(0,0,\dots,0)$ to $(0,1,0,\dots,0)$ in some $\varphi\colon \R^n\cong U\subseteq M$.
    
    If $x_0\neq x_1,x_2$, we can assume that $x_0\notin U$, we take $\sigma_1$ a stick avoiding $x_0,x_1,x_2$ in $M\setminus U$. If $x_0=x_1$ or $x_0=x_2$, we can first assume that $x_0=x_1=\varphi(0)$ and we fix $\sigma_1(t):=\varphi(t,0,\dots,0)$ (see Figure \ref{fig:path and stick}).
    \begin{figure}[H]
        \centering
        \includegraphics[width=0.1\linewidth]{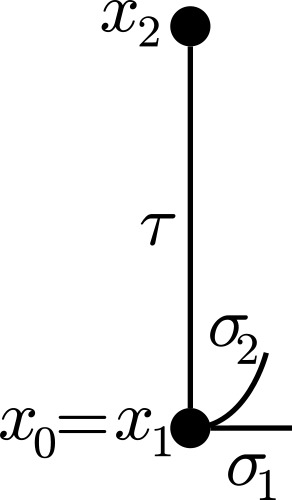}
        \caption{In the case $x_0=x_1$, the stick $\sigma_1$ is going to the right and the path $\tau$ and the $1$-parameter family is moving upwards in a local chart. Here we use the assumption that $n>1$.}
        \label{fig:path and stick}
    \end{figure}
    We then fix a $1$-parameter family supported in $U$ with $\phi^0=\id$ and that, in the ball $\varphi(B_2(0))$, is given by
    \begin{align*}
        \phi^s(\varphi(x_1,\dots,x_n))= \varphi(x_1,x_2+s,\dots, x_n).
    \end{align*}
    We then have that $\sigma_2(t):=\phi^{\chi(t)}(\sigma_1(t))$ is a stick avoiding $x_0$ and $x_2$. Therefore we can apply Proposition \ref{prop:a-s commuting with homotopy} and find 
    \begin{align*}
        (r_\tau\times l_{\overline{\tau}})\circ \lor^{x_1}_{\sigma_1}=\lor^{x_2}_{\sigma_2}.
    \end{align*}
    After dropping the dependence on $\sigma_1$ and $\sigma_2$, this gives the corollary.
\end{proof}

\begin{Cor}\label{cor:pi1 invariance}
    Assume that the dimension is $n>1$. The avoiding-stick coproduct $\lor$ factors through the fixed-points $H_*(\Omega M\times \Omega M)^{\pi_1(M,x_0)}$ where $\pi_1(M,x_0)$ acts on $\Omega M\times \Omega M$ via
    \begin{align*}
        [\tau]\cdot(\gamma,\gamma'):=(\gamma \cdot\tau, \tau^{-1} \cdot\gamma').
    \end{align*}
\end{Cor}
\begin{proof}
    This is the previous Corollary for the case $x_0=x_1=x_2$.
\end{proof}

\begin{Cor}\label{cor:infty pi1}
    Let $M$ be a manifold with infinite fundamental group of dimension $n>1$. The avoiding-stick coproduct and the Goresky-Hingston coproduct and the avoiding stick coproduct are zero on $H_*(\Omega M)$.
\end{Cor}
\begin{proof}
    By Proposition \ref{prop:a-s-coprod lifts GH}, it suffices to show that the avoiding-stick coproduct is zero.
    
    For all manifolds, it holds
    \begin{align*}
        H_*(\Omega M)&=\bigoplus_{\gamma\in \pi_1(M,x_0)}H_*(\Omega_\gamma M),\\
        H_*(\Omega M\times \Omega M )&=\bigoplus_{\gamma,\delta \in \pi_1(M,x_0)}H_*(\Omega_\gamma M\times \Omega_\delta M)
    \end{align*}
    where $\Omega_\gamma M$ is the connected component of $\Omega M$ of all loops that represent $\gamma\in \pi_1(M,x_0)$. The actions $r_\tau\times l_{\tau^{-1}}$ permutes the summands where $\gamma\delta$ are equal. Because $\pi_1(M,x_0)$ is infinite, there are infinite pairs in each orbit. Therefore the action has no fixed points except $0$. Corollary \ref{cor:pi1 invariance} now shows that the coproduct is zero.
\end{proof}

\begin{Prop}\label{prop:sullivan}
    Let $\sigma$ be a stick avoiding $x_0$. The Pontryagin product $\land$ and the avoiding-stick coproduct $\lor_\sigma$ satisfy the Sullivan relation:
    \begin{align*}
        \lor_\sigma\circ \land=(\id \times \land)\circ (\lor_\sigma\times \id)+(\land\times \id)\circ (\id \times \lor_\sigma).
    \end{align*}
\end{Prop}
\begin{proof}
    The strategy is to compute $\lor^{x}_{\widetilde{\sigma}}\circ\land(\alpha\times \beta)$ for a point $x\neq x_0$ and then use Proposition \ref{prop:a-s commuting with homotopy}.
    
    We fix a $1$-parameter family of homeomorphisms $\phi^s$ such that 
    \begin{align*}
        \widetilde{\sigma}\colon [0,\varepsilon]&\to M,\\
        s&\mapsto\phi^s\circ\sigma(s)
    \end{align*}
    is a stick avoiding $x_0$ and $x:=(\phi^1)^{-1}(x_0)\neq x_0$. Moreover, we fix a contractible neighbourhood $U_{x}$ such that $\widetilde{\sigma}(t)\notin U_{x}$ for all $t\in [0,\varepsilon]$. We set $\tau(t):=\phi^t(x)$. By Proposition \ref{prop:a-s commuting with homotopy}, we have
    \begin{align}\label{eq:tau and sigma for sull}
        (r_\tau\times l_{\overline{\tau}})\circ \lor_{\widetilde{\sigma}}^{x}=\lor_\sigma.
    \end{align}
    We define $\widetilde{\sigma}'(t):=\widetilde{\sigma}(2t)$ for $t\in [0,\frac{\varepsilon}{2}]$. We then have 
    \begin{align*}
        \land\circ (\pi_{\widetilde{\sigma}}\times \pi_{\widetilde{\sigma}})=\pi_{\widetilde{\sigma}'}\circ \land
    \end{align*}
    as maps $\Omega M\times \Omega M\to \Omega^{\widetilde{\sigma}'}M$. Moreover, the map
    \begin{align*}
        r':=r\circ (\id\times \land)\colon I \times \Omega^{\widetilde{\sigma}}M\times \Omega^{\widetilde{\sigma}}M\to \Omega M
    \end{align*}
    sends $([0,\frac{\varepsilon}{2}]\cup[1-\frac{\varepsilon}{2},1])\times \Omega^{\widetilde{\sigma}}M\times \Omega^{\widetilde{\sigma}}M$ to $\mathcal{F}_x^c$ because $\widetilde{\sigma}'$ avoids $x$. But because $\widetilde{\sigma}$ avoids $x$, $r'$ sends also $([\frac{1-\varepsilon}{2},\frac{1+\varepsilon}{2}])\times \Omega^{\widetilde{\sigma}}M\times \Omega^{\widetilde{\sigma}}M$ to $\mathcal{F}_x^c$.
    
    We therefore find that the inclusion
    \begin{align*}
        C_*\left(r'\left(((\tfrac{\varepsilon}{2},\tfrac{1-\varepsilon }{2})\cup(\tfrac{1+\varepsilon}{2},1-\tfrac{\varepsilon}{2}))\times \Omega^{\widetilde{\sigma}}M\times \Omega^{\widetilde{\sigma}}M\right),\mathcal{F}_x^c\right)\to C_*\left(r'\left(I\times \Omega^{\widetilde{\sigma}}M\times \Omega^{\widetilde{\sigma}}M\right),\mathcal{F}_x^c\right)
    \end{align*}
    is a homotopy equivalence by excision. We fix a homotopy inverse $\rho''$.
    
    Because $r'$ is the composition of two injective maps, it is injective and thus 
    \begin{align*}
        C_*&\left(r'\left(((\tfrac{\varepsilon}{2},\tfrac{1-\varepsilon }{2})\cup(\tfrac{1+\varepsilon}{2},1-\tfrac{\varepsilon}{2}))\times \Omega^{\widetilde{\sigma}}M\times \Omega^{\widetilde{\sigma}}M\right),\mathcal{F}_x\right)\\
        =&\ C_*\left(r'\left((\tfrac{\varepsilon}{2},\tfrac{1-\varepsilon }{2})\times \Omega^{\widetilde{\sigma}}M\times \Omega^{\widetilde{\sigma}}M\right),\mathcal{F}_x\right)\\
        &\oplus C_*\left(r'\left((\tfrac{1+\varepsilon}{2},1-\tfrac{\varepsilon }{2})\times \Omega^{\widetilde{\sigma}}M\times \Omega^{\widetilde{\sigma}}M\right),\mathcal{F}_x\right).
    \end{align*}
    Our computations so far show that $\lor^x_{\widetilde{\sigma}'}\circ \land$ is computed via the composition
    \begin{align*}
        C_*(\Omega M\times \Omega M)\overset{(\pi_{\widetilde{\sigma}}\times \pi_{\widetilde{\sigma}})_*}{\to }&\ C_*\left(\Omega^\sigma M\times \Omega^\sigma M\right)\overset{I\times}{\to} C_{1+*}(I\times \Omega^\sigma M\times \Omega^\sigma M, \partial I\times \Omega^\sigma M\times \Omega^\sigma M)\\
        \overset{r'_*}{\to}&\ C_{1+*}\left(r'\left(I\times \Omega^{\widetilde{\sigma}}M\times \Omega^{\widetilde{\sigma}}M\right),\mathcal{F}_x^c\right)\\
        \overset{\rho''}{\to}&\ C_{1+*}\left(r'\left((\tfrac{\varepsilon}{2},\tfrac{1-\varepsilon }{2})\times \Omega^{\widetilde{\sigma}}M\times \Omega^{\widetilde{\sigma}}M\right),\mathcal{F}_x^c\right)\\
        &\oplus C_{1+*}\left(r'\left((\tfrac{1+\varepsilon}{2},1-\tfrac{\varepsilon }{2})\times \Omega^{\widetilde{\sigma}}M\times \Omega^{\widetilde{\sigma}}M\right),\mathcal{F}_x^c\right)\\
        \overset{\lor_{\mathrm{triv}}^x}{\to}&\ C_{1-n+*}(\mathcal{P}_{x_0\to x}\times\mathcal{P}_{x\to x_0}).
    \end{align*}
    The injectivity of $r'$ also shows that for $\gamma_1, \gamma_2\in H_{*}( \Omega^\sigma M)$, we can represent $I\times \gamma_1\times \gamma_2\in H_{1+*}(I\times \Omega^\sigma M\times \Omega^\sigma M, \partial I\times \Omega^\sigma M\times \Omega^\sigma M)$ by
    \begin{align*}
        [\tfrac{\varepsilon}{2},\tfrac{1-\varepsilon}{2}]\times \gamma_1\times \gamma_2+[\tfrac{1+\varepsilon}{2},1-\tfrac{\varepsilon}{2}]\times \gamma_1\times \gamma_2+\alpha
    \end{align*}
    for some
    \begin{align*}
        \alpha&\in C_*((r')^{-1}(\mathcal{F}_x^c)).
    \end{align*}
    Because definition of $\lor^x_{\mathrm{triv}}$ is zero on $C_*(\mathcal{F}_x^c)$, we find $\lor^x_{\mathrm{triv}}(r'(\alpha))=0$.
    
    We also have that 
    \begin{align*}
        \lor^x_{\mathrm{triv}}(r'_*([\tfrac{\varepsilon}{2},\tfrac{1-\varepsilon}{2}]\times \gamma_1\times \gamma_2))=(\id\times\land)(\lor^x_{\widetilde{\sigma}}(\gamma_1)\times \gamma_2).
    \end{align*}
    Indeed after reparametrisation, we find
    \begin{align*}
        \lor^x_{\mathrm{triv}}(r'_*([\tfrac{\varepsilon}{2},\tfrac{1-\varepsilon}{2}]\times \gamma_1\times \gamma_2))=(\id\times\land)(\lor^x_{\mathrm{triv}}(r_*([{\varepsilon},{1-\varepsilon}]\times \gamma_1)\times \gamma_2))
    \end{align*}
    and $\lor_{\mathrm{triv}}^x(r_*([{\varepsilon},{1-\varepsilon}]\times \gamma_1))$ computes $\lor^x_{\widetilde{\sigma}}(\gamma_1)$.
    
    Similarly, we find
    \begin{align*}
        \lor^x_{\mathrm{triv}}(r'_*([\tfrac{1+\varepsilon}{2},1-\tfrac{\varepsilon}{2}]\times \gamma_1\times \gamma_2))=(-1)^{(n-1)|\gamma_1|}(\land\times \id)(\gamma_1\times\lor^x_{\widetilde{\sigma}}( \gamma_2))
    \end{align*}
    where the sign comes from the Koszul sign of moving the degree $-n$ map $\lor^x_{\mathrm{triv}}$ and the interval past $\gamma_1$.
    
    We thus have shown:
    \begin{align*}
        \lor^x_{\widetilde{\sigma}'}\circ \land(\gamma_1\times\gamma_2)=(\id\times\land)(\lor^x_{\widetilde{\sigma}}(\gamma_1)\times \gamma_2)+(-1)^{(n-1)|\gamma_1|}(\land\times \id)(\gamma_1\times\lor^x_{\widetilde{\sigma}}( \gamma_2)).
    \end{align*}
    After postcomposing with $r_\tau\times l_{\overline\tau}$, we find by \eqref{eq:tau and sigma for sull}
    \begin{align*}
        \lor_{\sigma}\circ \land(\gamma_1\times\gamma_2)=(\id\times\land)(\lor_{\sigma}(\gamma_1)\times \gamma_2)+(-1)^{(n-1)|\gamma_1|}(\land\times \id)(\gamma_1\times\lor_{\sigma}( \gamma_2)).
    \end{align*}
    This is the Sullivan relation evaluated on $\gamma_1\times \gamma_2$.
\end{proof}

\subsection{Higher Dimensional Spheres}
In this subsection, we consider the example of spheres of dimension $n>1$.

\begin{Exmp}[Spheres]\label{exmp:spheres}
    We recall that
    \begin{align*}
        H_*(\Omega S^n)\cong\Z[u]
    \end{align*}
    with $u$ in degree $n-1$. The right hand side is a free \textit{associative} algebras. In particular, the Pontryagin product is not (graded) commutative for $n$ even. The class $u$ corresponds to the class $au$ in \cite[Theorem 2 (2)]{cohen2004loop} and $AU$ in \cite[Subsection 8.2]{cieliebak2023loop} for $n$ odd.
    
    The class $u$ can be characterised as follows: we fix a vector $v_0\in UT_{x_0}S^n$. We have
    \begin{align*}
        D^{n-1}\cong K^{n-1}:= \{v\in U T_{x_0}\R^{n+1}\cong S^n\mid \langle  v_0,v\rangle=0, \langle v,x_0\rangle \leq 0\}.
    \end{align*}
    Here we interpret $U T_{x_0}\R^{n+1}\cong  S^{n}$ via the canonical embedding $S^{n}\to \R^{n+1}$. This is the space of all vectors pointing into the sphere which are perpendicular to $v_0$.
    
    A vector in $v\in K^{n-1}$ and $v_0$ define a plane $P_v$ containing $x_0$. There is (up to parametrization) a unique loop $\gamma_v\in  \Omega S^n$ in $P_v\cap S^{n}$ with start direction $v_0$. For $v\in \partial K^{n-1}$, we have that $P_v$ is tangential to $S^n$ and thus $P_v\cap S^n=x_0$. We then define $\gamma_v$ as the constant loop at $x_0$. See Figure \ref{fig:loop in sphere}.
    
    The class $u$ is represented by\footnote{The orientation on $S^n$ induces an orientation on $K^{n-1}/\partial K^{n-1}$: we have natural homeomorphisms $\langle v_0\rangle^\perp \cap S^n\cong S^{n-1}\cong K^{n-1}/\partial K^{n-1}$}:
    \begin{align*}
        S^{n-1}\cong K^{n-1}/\partial K^{n-1} &\to \Omega S^n,\\
        v& \mapsto \gamma_v.
    \end{align*}
    \begin{figure}[H]
        \centering
        \includegraphics[width=0.3\linewidth]{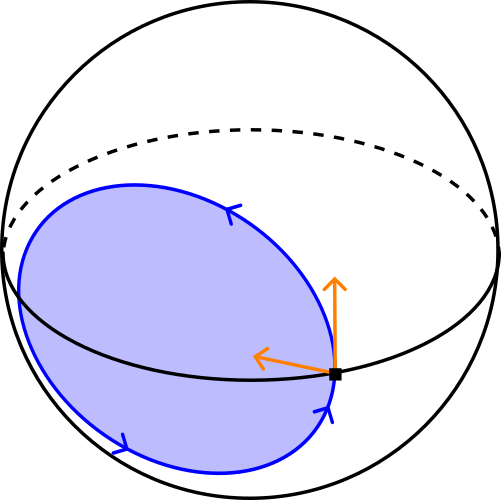}
        \caption{The vector $v_0$ pointing upward and a second perpendicular vector $v$ define a plane $P_v$ in bright blue. The plane $P_v$ intersects $S^n$ in one circle $\gamma_v$ in dark blue.}
        \label{fig:loop in sphere}
    \end{figure}
    For our computations, we fix $x_0=(1,0,\dots,0), x_1=(-1,\dots,0,0)\in S^n$. We choose the stick $\sigma(t)=(t,1-t^2,0,\dots,0)$ avoiding $x_0$ and $x_1$. Moreover, we fix $v_0=(0,\dots,0,1)\in UT_{x_0}S^n$. Then $K^{n-1}$ is given by
    \begin{align*}
        \{(v_0,\dots,v_{n-1},0)\in S^n\mid v_0\leq 0\}.
    \end{align*}
    For a $v=(v_0,\dots,v_{n-1},0)\in K^{n-1}$, we have
    \begin{align*}
        P_v:=\{(1,0,\dots, \lambda)+\mu (v_0,\dots,v_{n-1},0)\mid \lambda,\mu\in \R\}.
    \end{align*}
    In particular, solving 
    \begin{align*}
        x_1=(-1,\dots,0,0)=(1,0,\dots, \lambda)+\mu (v_0,\dots,v_{n-1},0)
    \end{align*}
    we find $\lambda=0$, $\mu=2$ and $v=(-1,0,\dots,0)$ as the only solution. In other words, the only $(t,v)\in I\times K^{n-1}$ such that $e(r(t,\sigma\overline\sigma\gamma_v\sigma\overline\sigma))=x_1$ are $t=\frac{1}{2}$ and $v=(-1,0,\dots,0)$.
    
    Moreover, this intersection is transversal because varying $t$ and varying $v$ gives the tangent space $T_{x_1}S^n$. Additionally, the local degree of these intersection is 1.
    
    Similarly, the only points $(t,v^{(1)},\dots,v^{(k)})\in I \times (K^{n-1})^k$ with $e(r(t,(\sigma\overline\sigma(\gamma_{v^{(1)}}\cdots\gamma_{v^{(k)}})\sigma\overline\sigma)))$ are 
    \begin{align*}
        t_i=\frac{i-1+\frac{1}{2}}{k}(1-4\varepsilon)+2\varepsilon,\text{ where } v^{(i)}=(-1,0,\dots,0)
    \end{align*}
    for $i=1,\dots,k$.
    
    We thus computed
    \begin{align*}
        \lor_{x_1}(u^k)=\sum_{i=1}^k u^{i-1} \cdot \gamma_{(-1,0,\dots,0)}|_{[0,\tfrac{1}{2}]}\otimes \gamma_{(-1,0,\dots,0)}|_{[\tfrac{1}{2},1]}\cdot u^{k-i}\in H_*(\mathcal{P}_{x_0\to x_1})\otimes H_*(\mathcal{P}_{x_1\to x_0}).
    \end{align*}
    We fix a path $\tau$ from $x_1$ to $x_0$. By Corollary \ref{cor:relate pre-coprod for n>1}, we find
    \begin{align*}
        \lor_{x_0}(u^k)=(r_{\tau}\otimes l_{\tau^{-1}})\circ \lor_{x_1}(u^k)=\sum_{i=1}^k u^{i-1} \cdot \gamma_{(-1,0,\dots,0)}|_{[0,\tfrac{1}{2}]}\cdot\tau\otimes \tau^{-1}\cdot\gamma_{(-1,0,\dots,0)}|_{[\tfrac{1}{2},1]}\cdot u^{k-i}.
    \end{align*}
    Because all loops in $S^n$ are contractible, so are the loops $\gamma_{(-1,0,\dots,0)}|_{[0,\tfrac{1}{2}]}\cdot\tau$ and $\tau^{-1}\cdot\gamma_{(-1,0,\dots,0)}|_{[\tfrac{1}{2},1]}$. Therefore we have
    \begin{align*}
        \lor_{x_0}(u^k)=\sum_{i=1}^k u^{i-1} \otimes u^{k-i}.
    \end{align*}
\end{Exmp}

\section{Maps and the Coproduct}\label{sec:Maps and Coprod}
In this section, we consider how our different coproducts interact with maps of manifolds. In Subsection \ref{subsec:univ covering}, we apply our results to the case of the universal covering map.\\

We start with the trivial coproduct and give a formula for it on chains. This lets us produce more interesting formulas for the avoiding-stick coproduct in Proposition \ref{prop:formula for f} and Corollaries \ref{cor:formula for f}, \ref{cor:deg 1}, \ref{cor:degree d}.

\begin{Lem}\label{lem:f and trivial coprod}
    Let $f\colon (M,x_0)\to (N,y_0)$ be a map of closed, oriented, pointed manifolds. Assume that $y_0$ is a regular value with $f^{-1}(y_0)=\{x_0,\dots, x_k\}$. There exist chain level coproducts $\lor_{\mathrm{triv},\Omega N}\colon C_*(\Omega N)\to C_*(\Omega N\times \Omega N)$ and $\lor^{x_i}_\mathrm{triv}\colon C_*(\Omega M)\to C_{-n+*}(\mathcal{P}_{x_0\to x_i}\times \mathcal{P}_{x_i\to x_0})$ for $i=0,\dots,k$ such that the following formula holds
    \begin{align*}
        \lor_{\mathrm{triv},\Omega N}\circ f_*=\sum_{i=0}^k\deg(f)_{x_i} (f\times f)_*\circ \lor^{x_i}_\mathrm{triv}
    \end{align*}
    as maps $C_*(\Omega M)\to C_{-n+*}(\Omega N\times \Omega N)$.
\end{Lem}
\begin{proof}
    Because $y_0$ is a regular value, we can find a contractible neighbourhood $y_0\in V\subseteq N$
    \begin{align*}
        f^{-1}(V)=\bigsqcup_{i=0}^k U_i
    \end{align*}
    where
    \begin{align*}
        f|_{U_i}\colon U_i\to V 
    \end{align*}
    is a homeomorphism and $x_i\in U_i$ for $i=0,\dots k$.
    
    We fix a trivial coproduct $\lor_{\mathrm{triv},\Omega N}\colon C_*(\Omega N)\to C_*(\Omega N\times \Omega N)$. This in particular fixes a homotopy retract $R_V\colon I\times V\to V$ fixing $y_0$, a local orientation $\tau_{y_0}\in C^n(V,\{y_0\}^c)$ and a homotopy equivalence
    \begin{align*}
        \rho^y\colon C_*(\Omega N,e^{-1}(y_0)^c)\to C_*(e^{-1}(V),e^{-1}(y_0)^c).
    \end{align*}
    We have that $R_{U_i}:=(\id\times f|_{U_i}^{-1})\circ R_V\circ (\id\times f)$ defines a homotopy retraction of $U_i$.
    
    We note that
    \begin{align*}
        f^{-1}(e^{-1}(V),e^{-1}(y_0)^c)=(e^{-1}(f^{-1}(V)),e^{-1}(f^{-1}(y_0))^c)=\bigsqcup_{i=0}^k (e^{-1}(U_i),e^{-1}(x_i)^c).
    \end{align*}
    We have
    \begin{align*}
        f^*\tau_{y_0}=\sum_{i=0}^k \deg(f)_{x_i}\tau_{x_i}
    \end{align*}
    for some local orientations $\tau_{x_i}\in C^n(U_i,\{x_i\}^c)$. This implies that
    \begin{align}\label{eq:pulling back local orien}
        f^*e^*\tau_{y_0}=\sum_{i=0}^k \deg(f)_{x_i}e^*\tau_{x_i}.
    \end{align}
    Moreover, we can find homotopy equivalences $\rho^{x_i}\colon C_*(\Omega M,e^{-1}(x_i)^c)\to C_*(e^{-1}(U_i),e^{-1}(x_i)^c)$ such that $f_*\circ \rho^{x_i}=\rho^{y_0}\circ f_*$.
    
    Therefore the following diagram commutes:
    \begin{center}
        \begin{tikzcd}
            C_*(\Omega M)\ar[r, "f_*"] \ar[d]& C_*(\Omega N)\ar[d]\\
            \bigoplus_{i=0}^k C_*(\Omega M,e^{-1}(x_i)^c)\ar[r, "f_*"]\ar[d,"\sum\rho^{x_i}"] & C_*(\Omega N,e^{-1}(y_0)^c)\ar[d, "\rho^{y_0}"]\\
            \bigoplus_{i=0}^k C_*(e^{-1}(U_i),e^{-1}(x_i)^c)\ar[r, "f_*"]\ar[d,"\sum\deg(f)_{x_i} e^*\tau_{x_i}\cap"] & C_*(e^{-1}(V),e^{-1}(y_0)^c)\ar[d, "e^*\tau_{y_0}\cap"]\\
            \bigoplus_{i=0}^k C_{-n+*}(e^{-1}(U_i))\ar[r, "f_*"]\ar[d,"R_*"] & C_{-n+*}(e^{-1}(V))\ar[d, "R_*"]\\
            \bigoplus_{i=0}^k C_{-n+*}(\mathcal{F}_{x_i}M)\ar[r, "f_*"]\ar[d,"c_*"] & C_{-n+*}(\mathcal{F}N)\ar[d, "c_*"]\\
            \bigoplus_{i=0}^k C_{-n+*}(\mathcal{P}_{x_0\to x_i}\times \mathcal{P}_{x_i\to x_0})\ar[r, "(f\times f)_*"] & C_{-n+*}(\Omega N\times \Omega N).
        \end{tikzcd}
    \end{center}
    The middle square commutes by naturality of the cap product and \eqref{eq:pulling back local orien}.
    
    This square proves the statement because the composition along the top right computes $\lor_{\mathrm{triv},\Omega N}\circ f_*$ while composition along the bottom left computes $\sum_{i=0}^k\deg(f)_{x_i} (f\times f)_*\circ \lor^{x_i}_\mathrm{triv}.$
\end{proof}
We next prove a formula for the avoiding-stick coproduct. This is a generalization of \cite[Lemma 3]{hingston2010loop}.
\begin{Prop}\label{prop:formula for f}
    Let $f\colon (M,x_0)\to (N,y_0)$ be a map of closed, oriented, pointed manifolds of dimension $n>1$. Assume that $y_0$ is a regular value with $f^{-1}(y_0)=\{x_0,\dots, x_k\}$. The following formula holds
    \begin{align*}
        \lor_{\Omega N}\circ f_*=\sum_{i=0}^k\deg(f)_{x_i} (f\times f)_*\circ \lor^{x_i}
    \end{align*}
    as maps $H_*(\Omega M)\to H_{1-n+*}(\Omega N\times \Omega N)$.
\end{Prop}
\begin{proof}
    We can fix $\sigma$ a stick avoiding $x_0,\dots, x_k$ such that $f\circ \sigma$ is a stick avoiding $y_0$.
    
    We have that $f\circ \pi_\sigma=\pi_{f\circ \sigma}\circ f$ as maps $\Omega M\to \Omega^{f\circ \sigma}N$. Moreover, $f$ commutes with the map $I\times $ and intertwines the maps
    \begin{align*}
        r\colon (I\times \Omega^\sigma M,\partial I\times \Omega^\sigma M)\to (r(I\times \Omega^\sigma M),e^{-1}(\{x_0,\dots,x_k\})^c\cup \mathcal{H})
    \end{align*}
    and 
    \begin{align*}
        r\colon (I\times\Omega^{f\circ \sigma}N,\partial I\times \Omega^{f\circ \sigma} N)\to (r(I\times \Omega^{f\circ \sigma}N),e^{-1}(y_0)^c\cup \mathcal{H}).
    \end{align*}
    Similarly, $f_*$ intertwines the map
    \begin{align*}
        C_*(r((\varepsilon,1-\varepsilon)\times \Omega^\sigma M),e^{-1}(\{x_0,\dots,x_k\})^c)\to C_*(r(I\times \Omega^\sigma M),e^{-1}(\{x_0,\dots,x_k\})^c\cup \mathcal{H})
    \end{align*}
    and 
    \begin{align*}
        C_*(r((\varepsilon,1-\varepsilon)\times \Omega^{f\circ \sigma}N),e^{-1}(y_0)^c)\to C_*(r(I\times \Omega^{f\circ \sigma} M),e^{-1}(y_0)^c\cup \mathcal{H}).
    \end{align*}
    Therefore their inverses in homology $\sum_i \rho_\sigma^{x_i}$ and $\rho_{f\circ \sigma}$ are intertwined by $f_*$ in homology.
    
    We thus have shown that
    \begin{align*}
        \lor_{f\circ \sigma}\circ f_*&=\lor_{\mathrm{triv}}\circ \rho_{f\circ\sigma}\circ r_*\circ (I\times)\circ (\pi_{f\circ \sigma})_* \circ f_*\\
        &=\lor_{\mathrm{triv}}\circ f_*\circ \left(\sum_{i=0}^k \rho_{\sigma}^{x_i}\right)\circ r_*\circ (I\times)\circ (\pi_{\sigma})_*.
    \end{align*}
    Thus by Lemma \ref{lem:f and trivial coprod}, we find
    \begin{align*}
        \lor_{f\circ \sigma}\circ f_*&=\left(\sum_{i=0}^k\deg(f)_{x_i} (f\times f)_*\circ \lor^{x_i}_\mathrm{triv}\right)\circ  \left(\sum_{i=0}^k \rho_{\sigma}^{x_i}\right)\circ r_*\circ (I\times)\circ (\pi_{\sigma})_*.
    \end{align*}
    Finally, we observe that, if $i\neq j$, $\lor_{\mathrm{triv}}^{x_i}$ is zero on $C_*(r((\varepsilon,1-\varepsilon)\times \Omega^\sigma M),e^{-1}(x_j)^c)$  and thus $\lor_{\mathrm{triv}}^{x_i}\circ \rho_\sigma^{x_j}=0$. This lets us compute:
    \begin{align*}
        \lor_{f\circ \sigma}\circ f_*&=\left(\sum_{i=0}^k\deg(f)_{x_i} (f\times f)_*\circ \lor^{x_i}_\mathrm{triv}\circ \rho_{\sigma}^{x_i}\right)\circ r_*\circ (I\times)\circ (\pi_{\sigma})_*\\
        &=\sum_{i=0}^k\deg(f)_{x_i} (f\times f)_*\circ \lor^{x_i}_\sigma.
    \end{align*}
    Because we work in homology and dimension $n>1$, we can drop the dependence on the stick which gives the desired formula (see Corollary \ref{cor:independence of sigma}).
\end{proof}

Using the fact that we can relate pre-coproducts to the usual coproduct, we can find the following clean formula.

\begin{Cor}\label{cor:formula for f}
    Let $f\colon (M,x_0)\to (N,y_0)$ be a map of closed, oriented, pointed manifolds of dimension $n>1$. Assume that $y_0$ is a regular value with $f^{-1}(y_0)=\{x_0,\dots, x_k\}$. Let $\tau_i$ be paths from $x_i$ to $x_0$. For the avoiding-stick coproducts, it holds
    \begin{align*}
        \lor_{\Omega N} \circ f_*=\sum_{i=0}^k \deg(f)_{x_i} (r_{[f_*(\tau_i)]}\times l_{[f_*(\tau_i)]^{-1}})\circ (f\times f)_*\circ \lor_{\Omega M}
    \end{align*}
    as maps $H_*(\Omega M)\to H_{1-n+*}(\Omega N\times \Omega N)$
\end{Cor}
\begin{proof}
    By Proposition \ref{prop:formula for f}, we have
    \begin{align*}
        \lor_{\Omega N}\circ f_*=\sum_{i=0}^k\deg(f)_{x_i} (f\times f)_*\circ \lor^{x_i}.
    \end{align*}
    By Corollary \ref{cor:relate pre-coprod for n>1}, we have
    \begin{align*}
        \lor^{x_i}=(r_{\tau_i}\times l_{\overline{\tau_i}})\circ \lor^{x_0}=(r_{\tau_i}\times l_{\overline{\tau_i}})\circ \lor_{\Omega M}.
    \end{align*}
    We further note that
    \begin{align*}
        (f_*\times f_*)\circ (r_{\tau_i}\times l_{\overline{\tau_i}})=(r_{f\circ \tau_i}\times l_{\overline{f\circ \tau_i}})\circ (f_*\times f_*).
    \end{align*}
    Here $f\circ\tau_i$ is a loop and it thus represents an element in $\pi_1(N,y_0)$. In this group, it holds $[\overline{f\circ \tau_i}]=[{f\circ \tau_i}]^{-1}$. Therefore, we have
    \begin{align*}
        \lor_{\Omega N}\circ f_*&=\sum_{i=0}^k\deg(f)_{x_i} (f\times f)_*\circ \lor^{x_i}\\
        &=\sum_{i=0}^k\deg(f)_{x_i} (f_*\times f_*)\circ (r_{\tau_i}\times l_{\overline{\tau_i}})\circ \lor_{\Omega M}\\
        &=\sum_{i=0}^k\deg(f)_{x_i} (r_{[f\circ \tau_i]}\times l_{[f\circ \tau_i]^{-1}})\circ (f_*\times f_*)\circ \lor_{\Omega M}.
    \end{align*}
\end{proof}

In particular, for degree $1$ maps we get that the Goresky-Hingston coproduct commutes with them:

\begin{Cor}\label{cor:deg 1}
    Let $f\colon (M,x_0)\to (N,y_0)$ be a degree 1 map of closed, oriented, pointed manifolds which preserves orientation. Then $(f\times f)_*\circ \lor_{\Omega M}=\lor_{\Omega N}\circ f_*$ holds.
\end{Cor}
\begin{proof}
    If $n=1$, then $f$ is homotopic to a homeomorphism and we are done.
    
    If $n>1$, we can replace $f$ with a homotopic map $g$ such that $g^{-1}(y_0)=\{x_0\}$ (see \cite{hopf1930topologie} for the case $n>2$ and \cite{kneser1928glattung} for $n=2$). We thus have $\deg(g)_{x_0}=1$. We then choose $\tau_0$ the constant path at $x_0$ and apply Corollary \ref{cor:formula for f} which concludes the proof.
\end{proof}

\begin{Cor}\label{cor:degree d}
    Let $(M,x_0)$ and $(N,y_0)$ be closed, oriented, pointed manifolds of dimension $n>1$ and $N$ simply connected. Let $f\colon (M,x_0)\to (N,y_0)$ be a map of degree $d$. It holds $\lor_{\Omega N}\circ f_*=d\cdot (f\times f)_*\circ \lor_{\Omega M}$.
\end{Cor}
\begin{proof}
    Because $N$ is simply connected, we have that $f_*(\tau_i)\in\pi_1(N,y_0)$ is trivial for all choices $\tau_i$. We apply Corollary \ref{cor:formula for f} and find
    \begin{align*}
        \lor_{\Omega N}\circ f_*=\left(\sum_{i=0}^k\deg(f)_{x_i} \right)(f\times f)_*\circ \lor_{\Omega N}=d\cdot (f\times f)_*\circ \lor_{\Omega N}.
    \end{align*}
\end{proof}

\subsection{Universal Covering}\label{subsec:univ covering}
In this subsection, we study the case where the map $p\colon \widetilde{M}\to M$ we consider is the universal covering. In this special case, we can describe the coproduct on $M$ in terms of the coproduct of $\widetilde{M}$ and vice-versa.\\

We first recall that, as chain complexes, the chains on the based loop space of a space is described by the chains of the based loop space of the universal covering.

\begin{Lem}\label{lem:Omega X with univ cov}
    Let $(X,x_0)$ be a pointed, connected CW complex with universal covering $p\colon (\widetilde{X},\widetilde{x}_0)\to (X,x_0)$. The chains of the based loop space of $X$ has the following description:
    \begin{align*}
        C_*(\Omega_{x_0}{X})\simeq \bigoplus_{l\in \pi(X,x_0)} C_*(\Omega_{\widetilde{x}_0} \widetilde{X}).
    \end{align*}
\end{Lem}
\begin{proof}
    We know that $\Omega X$ has connected components corresponding to $l\in \pi_1 (X)$:
    \begin{align*}
        \Omega_l X:=\{\gamma\in \Omega X\mid [\gamma]=l\in\pi_1 (X)\}.
    \end{align*}
    Because $\Omega X$ is group-like, we have that $\Omega_l X\simeq \Omega_{1}X$ for $1\in \pi_1 (X)$ the neutral element.
    
    On the other hand, we know that the map $p_*\colon \pi_k(\widetilde X)\to \pi_k (X)$ is an isomorphism for $k\geq 2$. Thus for $k\geq 1$, the map that $p$ induces on $\pi_k(\Omega\widetilde X)$ is given by the isomorphism
    \begin{align*}
        \pi_k(\Omega \widetilde X,\widetilde{x}_0)\cong \pi_{k+1} (\widetilde X,\widetilde{x}_0)\overset{\pi_{k+1}(p)}{\cong}\pi_{k+1}(X,x_0)\cong \pi_k(\Omega X,x_0)\cong \pi_k(\Omega_1X).
    \end{align*}
    Therefore $p$ induces a homotopy equivalence $\Omega \widetilde{X}\cong \Omega_1 X$. This concludes the proof.
\end{proof}

\begin{Rem}\label{rem:Omega X univ cov}
    The proof of Lemma \ref{lem:Omega X with univ cov} gives a description of $H_*(\Omega X)$ in the above setting: any element $\alpha\in H_*(\Omega X)$ can be uniquely written as
    \begin{align*}
        \alpha=\sum_{g\in \pi_1 X} g\cdot p_*(\alpha_g)
    \end{align*}
    for some $\alpha_g\in H_*(\Omega \widetilde{X})$ only finitely many non-zero.
\end{Rem}

We now turn our attention back to $f\colon \widetilde{M}\to M$, the universal covering of a closed manifold. We want to understand the Goresky-Hingston coproduct on $\Omega M$ in terms of the coproduct on $\Omega \widetilde{M}$. We recall that by Corollary \ref{cor:infty pi1} the coproduct is trivial if $\pi_1(M,x_0)$ is infinite. We can thus assume that $\pi_1(M,x_0)$ is finite. Therefore the covering space $\widetilde{M}$ is itself a closed manifold and has indeed a Goresky-Hingston coproduct on its based loop space. We get the following description of the coproduct on $\Omega M$:
\begin{Prop}\label{prop:GH via univ cov}
    Let $(M,x_0)$ be a pointed, closed, oriented manifold with finite fundamental group and universal covering $p\colon (\widetilde{M},y_0)\to (M,x_0)$. Let $\alpha=g\cdot p_*(\beta)$ for some $g\in \pi_1(M,x_0)$ and $\beta\in H_*(\Omega \widetilde M)$ with avoiding-stick coproduct:
    \begin{align*}
        \lor_{\Omega\widetilde{M}}(\beta)=\sum_i\gamma_{i}\times \delta_{i}
    \end{align*}
    The avoiding-stick coproduct on $\alpha$ is computed as
    \begin{align*}
        \lor_{\Omega M}(\alpha)=\sum_{h\in \pi_1(M,x_0)} \sum_i  g\cdot p_*(\gamma_{i})\cdot h\times h^{-1}\cdot p_*(\delta_{i}).
    \end{align*}
\end{Prop}
\begin{proof}
    By the Sullivan relation as in Proposition \ref{prop:sullivan}, we have
    \begin{align*}
        \lor_{\Omega M}(g\cdot p_*(\beta))=\lor_{\Omega M}(g)\cdot p_*(\beta)+g\cdot \lor_{\Omega M}(p_*(\beta))=0+g\cdot \lor_{\Omega M}(p_*(\beta)).
    \end{align*}
    We now apply the formula from Corollary \ref{cor:formula for f} to compute $\lor_{\Omega M}( p_*(\beta))$:
    \begin{align*}
        \lor_{\Omega M} \circ p_*=\sum_{x_i\in p^{-1}(y_0)} \deg(p)_{x_i} (r_{[p_*(\tau_i)]}\times l_{[p_*(\tau_i)]^{-1}})\circ (p \times p)_*\circ \lor_{\Omega \widetilde{M}},
    \end{align*}
    where $\tau_i$ are paths connecting $x_0$ and $x_i$. We note that in our case all $x_i$ have local degree $+1$ and the loops $[p_*(\tau_i)]$ correspond to all elements $h\in \pi_1(M)$. Therefore, we find
    \begin{align*}
        \lor_{\Omega M}( p_*(\beta))=\sum_{h\in \pi_1(M)} \sum_i  p_*(\gamma_{i})\cdot h\times h^{-1}\cdot p_*(\delta_{i})
    \end{align*}
    which concludes the proof.
\end{proof}

\section{Spheres Modulo a Group Action}\label{sec:Sn/G}
In this section, we apply our methods to $M=S^n/G$ for a finite group $G$. We can fully describe the Goresky-Hingston coproduct on $\Omega M$. This also helps us describe the homology of the free loop space because the coproduct gives non-vanishing results on the relevant Leray-Serre spectral sequence.\\ 

Let $G$ be a non-trivial, finite group that acts freely and orientation-preserving on $S^n$ for some $n\geq 2$. We note that $p\colon S^n\to S^n/G$ is the universal covering of the oriented manifold $M:=S^n/G$.

We have the following description of the Pontryagin ring:
\begin{Prop}\label{prop:Omega Sn/G}
    Let $G$ be a non-trivial, finite group that acts on $S^n$ freely and orientation-preserving for some $n\geq 2$ and let $M=S^n/G$. The Pontryagin ring on $\Omega M$ is given by
     \begin{align*}
         H_*(\Omega M)\cong \Z[G][x]
     \end{align*}
     where $\Z[G]$ is the group ring and $x$ is a central element of degree $n-1$.
\end{Prop}
\begin{proof}
    By Lemma \ref{lem:Omega X with univ cov}, we find that as abelian groups
    \begin{align*}
        H_*(\Omega M)\cong \bigoplus_G H_*(S^n).
    \end{align*}
    We recall that $H_*(S^n)\cong \Z[u]$ for $u$ a class of dimension $n-1$. We thus set $x:=p_*(u)\in H_{n-1}(\Omega_1 M)$ and find
    \begin{align*}
        H_*(\Omega M)\cong \Z[G][x]
    \end{align*}
    as abelian groups.
    
    Because $n$ is odd, $u\in H_*(S^n)$ commutes with itself and therefore so does $x$. It thus remains to show that $x$ and $g\in G$ commute. Because conjugation by $g$ is a $G$-action on $H_{n-1}(\Omega_1 M)\cong \Z$, we have $g^{-1} xg=\pm x$.
    
    We thus have to show that the sign is $+1$. Morally this is due to the fact that the class $u$ comes from the orientation of $S^n$ and $g$ acts trivially on the orientation. We spell this argument out in more detail.
    
    We have that
    \begin{center}
        \begin{tikzcd}
            H_{n-1}(\Omega_{x_0}S^n)\cong H_n(\Sigma \Omega_{x_0}S^n) \ar[r, "\mathrm{ev}"]& H_n(S^n)
        \end{tikzcd}
    \end{center}
    sends $u$ to the orientation $[S^n]\in H_n(S^n)$.
    
    The diagram 
    \begin{center}
        \begin{tikzcd}
            \Sigma\Omega_{x_0} S^n \ar[r, "\mathrm{ev}"] \ar[d,"p"] & S^n \ar[d, "p"]\\
            \Sigma \Omega M \ar[r, "\mathrm{ev}"] & M
        \end{tikzcd}
    \end{center}
    commutes. Because we have $p_*(S^n)=|G|\cdot [M]$, we find $\mathrm{ev}_*(x)=|G|\cdot [M]$. It thus suffice to show that $\mathrm{ev}_*(g^{-1}xg)=|G|\cdot [M]$.
    
    The commuting diagram
    \begin{center}
        \begin{tikzcd}
            \Sigma \Omega_{x_0} S^n \ar[r] \ar[d, "g"] & S^n \ar[d,"g"]\\
            \Sigma \Omega_{g\cdot  x_0}S^n \ar[r] & S^n
        \end{tikzcd}
    \end{center}
    shows that $\mathrm{ev}_*$ sends $g_*( u)\in H_{n-1}(\Omega_{g \cdot x_0}S^n)$ to $(\varphi_g)_*([S^n])=[S^n]$.
    
    We now give a different description of $g_*(u)$. Let $\alpha\in C_{n-1}(\Omega_{x_0}S^n)$ be a representative of $u$ and let $\tau_g$ be a path from $x_0$ to $g\cdot x_0$. Then $\tau_{g}^{-1}\cdot\alpha\cdot \tau_g$ represents $g_*( u)$. We thus have
    \begin{align*}
        p_*( g_*(u))=p_*(\tau_{g}^{-1})\cdot p_*([\alpha])\cdot p_*(\tau_g)=|G|g^{-1}\cdot x\cdot g.
    \end{align*}
    Now the commuting diagram
    \begin{center}
        \begin{tikzcd}
            \Sigma\Omega_{g\cdot x_0} S^n \ar[r, "\mathrm{ev}"] \ar[d,"p"] & S^n \ar[d, "p"]\\
            \Sigma \Omega M \ar[r, "\mathrm{ev}"] & M
        \end{tikzcd}
    \end{center}
    shows that 
    \begin{align*}
        \mathrm{ev}_*(g^{-1}ug)=\mathrm{ev}_*(p_*(g_*( u)))=p_*(\mathrm{ev}_*(g_*( u)))=p_*([S^n])=|G|\cdot[M].
    \end{align*}
    This in particular shows that $g^{-1}xg=x$ and thus concludes the proof.
\end{proof}

\begin{Cor}\label{cor:coprod Omega Sn/G}
    Let $G$ be a non-trivial, finite group that acts on $S^n$ freely and orientation-preserving for some $n\geq 2$ and let $M=S^n/G$. Let $gx^k\in H_*(\Omega M)$. The avoiding-stick coproduct extends to a map $H_{n-1+*}(\Omega M)\to H_*(\Omega M)\otimes H_*(\Omega M)$ which is computed as
    \begin{align*}
        \lor(gx^k)=\sum_{i+j=k-1}\sum_{h\in G} gh^{-1}x^i\otimes hx^j.
    \end{align*}
\end{Cor}
\begin{proof}
    Because all homology groups are free, the Künneth map $H_*(\Omega M)\otimes H_*(\Omega M)\to H_*(\Omega M\times \Omega M)$ is an isomorphism and we can thus postcompose the avoiding-stick coproduct $H_{n-1+*}(\Omega M)\to H_*(\Omega M\times \Omega M)$ with the inverse of the Künneth map and get a map $H_{n-1+*}(\Omega M)\to H_*(\Omega M)\otimes H_*(\Omega M)$.
    
    The formula then follows from the coproduct on $H_*(\Omega S^n)$ as described in Example \ref{exmp:spheres} transported to $H_*(\Omega M)$ using Proposition \ref{prop:GH via univ cov} and Proposition \ref{prop:Omega Sn/G}.
\end{proof}

\subsection{The Free Loop Space}
We now turn our attention to the free loop space $\Lambda M$. We show that the Leray-Serre spectral sequence associated to the fibration $\mathrm{ev}_0\colon \Lambda M\to M$ is such that the $E^2$-page corresponds to the $E^\infty$-page. We do this using the description of the coproduct on $H_*(\Omega M)$ as in Corollary \ref{cor:coprod Omega Sn/G}.\\

We first compute $H_*(M)$: we consider the fibration
\begin{center}
    \begin{tikzcd}
        S^n \ar[r] & M\ar[d]\\
        & BG.
    \end{tikzcd}
\end{center}
The associated Leray-Serre spectral sequence collapses to the long exact sequence:
\begin{center}
    \begin{tikzcd}
        \dots \ar[r] & H_{i+1}(G) \ar[r] & H_{i+1-n}(G) \ar[r]& H_i(M) \ar[r] &H_{i}(G) \ar[r] & H_{i-n}(G) \ar[r]&\dots
    \end{tikzcd}
\end{center}
We know that $H_0(M)=H_n(M)=\Z$. The above long exact sequence then additionally implies $H_i(M)\cong H_i(G)$ for $1\leq i\leq n-1$.\\

Now we can study the space $\Lambda M$. It has connected components corresponding to conjugacy classes in $G$. For $[\gamma]\in \mathrm{ccl}(G)$, we denote $\Lambda_{[g]}M\subseteq \Lambda M$ for the corresponding connected component. It has associated the following fibration:
\begin{equation}\label{eq:spec seq Lambda[g]M}
    \begin{tikzcd}
        \bigsqcup_{h\in [g]} \Omega_h M \ar[r]& \Lambda_{[g]}M\ar[d, "\mathrm{ev}_0"]\\
        & M.
    \end{tikzcd}
\end{equation}
We note that as $G$-modules, we have
\begin{align*}
    H_*\left(\bigsqcup_{h\in [g]} \Omega_h M\right)\cong \bigoplus_{h\in [g]}H_*(\Omega_hM)\cong \begin{cases}
        \Z[G/C_G(g)] &\text{if }*=m(n-1),\, m\geq 0;\\
        0 &\text{else.}
    \end{cases}
\end{align*}
Therefore, the $E^2$-page of the associated Leray-Serre spectral sequence thus reads as 
\begin{align*}
    E^2_{p,q}([g]):=\begin{cases}
        H_p(M;\Z[G/C_G(g)]) &\text{if } q=m(n-1),\, m\geq 0;\\
        0& \text{else.}
    \end{cases}
\end{align*}
We denote the corresponding higher pages of the spectral sequence by $E^r_{p,q}([g])$.

By \cite[Example 3H.2]{hatcher2002algebraic} and our previous computations, we find
\begin{align*}
    H_p(M;\Z[G/C_G(g)])\cong H_p(S^n/C_G(g))\cong \begin{cases}
        \Z & \text{if } p=0,n;\\
        H_p(C_G(g))&\text{if }1\leq p\leq n-1;\\
        0&\text{else.}
    \end{cases}
\end{align*}
We thus have found that the $E^2$-page of Leray-Serre of $E^r_{p,q}([g])$ looks like this.
\begin{center}
    \begin{tikzpicture}
    \matrix (m) [matrix of math nodes, nodes in empty cells, row sep=1ex, column sep=1ex, minimum width=10ex]
        {
            \vdots\hspace{-1.5cm} &\vdots &   \vdots &  \vdots & \vdots & \vdots&\vdots\\ 
            n-1 \hspace{-1.5cm}&\Z& H_1(C_G(g)) &  \cdots & H_{n-1}(C_G(g))  & \Z &0\\ 
            n-2\hspace{-1.5cm}&0 &  0&   0  & 0& 0 &0 \\
            \vdots\hspace{-1.5cm} &\vdots &   \vdots &  \vdots & \vdots & \vdots&\vdots\\ 
            1\hspace{-1.5cm}&0 &  0&   0  & 0& 0 &0 \\
            0\hspace{-1.5cm}&\Z & H_1(C_G(g))  & \cdots & H_{n-1}(C_G(g))  & \Z &0 \\ 
            &0 &  1  & \dots &n-1 &n &>n \\ 
        };
        \path[-stealth] 
        ;
        \draw (-4.25,-2.5) -- (-4.25,2.5); 
        \draw (-5.35,-2) -- (6.1,-2); 
    \end{tikzpicture}
\end{center}

\begin{Prop}\label{prop:E2 is Einfty}
    For all $g\in G$ it holds that $E^2_{*,*}([g])\cong E^\infty_{*,*}([g])$.
\end{Prop}
\begin{proof}
    For degree reasons, the only differential that may show up is $d^n\colon E^n_{n,2k}([g])\to E^n_{0,(n-1)k}([g])$. Thus it suffices to show that $E^\infty_{0,(n-1)k}([g])\cong  E^2_{0,(n-1)k}([g])$ for $k>0$. By our computations, we have $E^2_{0,(n-1)k}([g])\cong \Z$. Moreover because all other groups of total degree of $(n-1)k$ are torsion, it suffices to show that $\iota_*\colon H_{(n-1)k}(\Omega_gM)\to H_{(n-1)k}(\Lambda_{[g]} M,M)$ is injective for $g\neq 1$ or $k>0$.
    
    We do this by induction on $k$. The statement holds for $k=0$ and $g\neq 1$ because $H_0(\Lambda_{[g]}M)\cong\Z$.
    
    For the induction step, we use the Goresky-Hingston coproduct on $H_*(\Lambda M,M)$ and the avoiding-stick coproduct on $H_*(\Omega M)$. The avoiding-stick coproduct is a lift of the Goresky-Hingston coproduct on the based loop space, which is a restriction of the Goresky-Hingston coproduct on the free loop space. We thus have the commuting diagram
    \begin{center}
        \begin{tikzcd}
            H_{k(n-1)}(\Omega_g M) \ar[r, "{\iota_*}"] \ar[d, "{\lor_{\Omega}}"] & H_{k(n-1)}(\Lambda_{[g]} M,M) \ar[d, "\lor_{\mathrm{GH}}"]\\
            H_{(k-1)(n-1)}(\Omega M\times \Omega M) \ar[r, "{(\iota\times \iota)_*}"] & H_{(k-1)(n-1)}((\Lambda M,M)\times (\Lambda M,M)).
        \end{tikzcd}
    \end{center}
    Our strategy is to find an element in $H_{k(n-1)}(\Omega_g M)$ that gets send to a non-torsion element in $ H_{(k-1)(n-1)}((\Lambda M,M)\times (\Lambda M,M))$ under composition on the bottom left. This then shows that there exists a non-torsion element in the image of the top map which is what we want to show.
    
    By induction, we can assume that, for $k_1+k_2=k-1$, the map
    \begin{align*}
        H_{k_1(n-1)}(\Omega_{g_1} M)\otimes H_{k_2(n-1)}(\Omega_{g_2} M)\to H_{(k-1)(n-1)}((\Lambda M,M)\times (\Lambda M,M))
    \end{align*}
    is injective if $k_j\geq 1$ or $g_j\neq 1$ for $j=1,2$.
    
    Corollary \ref{cor:coprod Omega Sn/G} then gives on $gx^k\in H_{k(n-1)}(\Omega_gM)$
    \begin{align*}
        (\iota\times \iota)_*\circ\lor_{\Omega M}(gx^k)=\sum_{g_1\cdot g_2=g} \sum_{k_1+k_2=k-1} \iota_*(g_1x^{k_1-1})\times \iota_*(g_2x^{k_2-1})
    \end{align*}
    which is non-zero if there exists $g_1,g_2\neq 1$ with $g_1g_2=g$ or $k_1,k_2>1$. This is satisfied except in the case where $G=\{\pm 1\}$, $g=-1$ and $k=1$.
    
    The remaining case is for the manifold $M=\R P^n$. We thus need to show that $H_{n-1}(\Lambda_{[-1]} \R P^n,\R ^n)$ has a non-torsion summand. This can be seen from the fact that the Chas-Sullivan product with $x$ gives an isomorphism $H_{*}(\Lambda \R P^n,\R P^n)\cong H_{n-1+*}(\Lambda \R P^n,\R P^n)$ \cite[1.11]{goresky2009loop} and our computation that shows that  $H_{2(n-1)}(\Lambda_{[-1]} \R P^n,\R ^n)$ has a non-torsion summand.
\end{proof}

\begin{Cor}\label{cor:free loop space}
    The homology of the free loop space of $M=S^n/G$ fits in the following short exact sequence:
    \begin{center}
        \begin{tikzcd}
            0 \ar[r] & H_{i-{n-1}}(S^n/C_G(g)) \ar[r]& H_{i+k(n-1)}(\Lambda_{[g]} M) \ar[r] & H_i(S^n/C_G(g)) \ar[r] &0
        \end{tikzcd}
    \end{center}
    for $k\geq 0$ and $1\leq i \leq n$. This short exact sequence splits for all $i\neq n-1$.
\end{Cor}
\begin{proof}
    The short exact sequence comes immediately from the description $E^\infty_{*,*}([g])$ as in Proposition \ref{prop:E2 is Einfty}. Moreover, for $i< n-1$ we have $H_{i-{n-1}}(S^n/C_G(g))\cong 0$. For $i=n$, we have $H_i(S^n/C_G(g))\cong \Z$ is projective and thus the short exact sequence splits.
\end{proof}

\normalem
\bibliographystyle{alpha}
\bibliography{ref.bib}

\begin{thebibliography}{CHO23}

\bibitem[CHO23]{cieliebak2023loop}
Kai Cieliebak, Nancy Hingston, and Alexandru Oancea.
\newblock Loop coproduct in {M}orse and {F}loer homology.
\newblock {\em Journal of Fixed Point Theory and Applications}, 25(2):59, 2023.

\bibitem[CJY04]{cohen2004loop}
Ralph~L Cohen, John~DS Jones, and Jun Yan.
\newblock The loop homology algebra of spheres and projective spaces.
\newblock In {\em Categorical Decomposition Techniques in Algebraic Topology: International Conference in Algebraic Topology, Isle of Skye, Scotland, June 2001}, pages 77--92. Springer, 2004.

\bibitem[Cli25]{clivio2025GHviaDGMorse}
Jonathan Clivio.
\newblock {The Goresky-Hingston Coproduct in Morse Homology with DG Coefficients}.
\newblock {\em In preparation}, 2025.

\bibitem[EK71]{edwards1971deformations}
Robert~D Edwards and Robion~C Kirby.
\newblock Deformations of spaces of imbeddings.
\newblock {\em Annals of Mathematics}, 93(1):63--88, 1971.

\bibitem[GH09]{goresky2009loop}
Mark Goresky and Nancy Hingston.
\newblock Loop products and closed geodesics.
\newblock {\em Duke Mathematical Journal}, 150(1):117, 2009.

\bibitem[Hat02]{hatcher2002algebraic}
Allen Hatcher.
\newblock {\em {Algebraic Topology}}.
\newblock Cambridge University Press, 2002.

\bibitem[Hin10]{hingston2010loop}
Nancy Hingston.
\newblock Loop products on connected sums of projective spaces.
\newblock {\em A celebration of the mathematical legacy of Raoul Bott}, 50:161--175, 2010.

\bibitem[Hop30]{hopf1930topologie}
Heinz Hopf.
\newblock {Zur Topologie der Abbildungen von Mannigfaltigkeiten: Zweiter Teil Klasseninvarianten von Abbildungen}.
\newblock {\em Mathematische Annalen}, 102(1):562--623, 1930.

\bibitem[HW17]{hingston2017product}
Nancy Hingston and Nathalie Wahl.
\newblock {Product and coproduct in string topology}.
\newblock {\em arXiv preprint arXiv:1709.06839}, 2017.

\bibitem[Kne28]{kneser1928glattung}
Hellmuth Kneser.
\newblock {Gl{\"a}ttung von Fl{\"a}chenabbildungen}.
\newblock {\em Mathematische Annalen}, 100(1):609--617, 1928.

\bibitem[KP24]{kenigsberg2024obstructions}
Lea Kenigsberg and Noah Porcelli.
\newblock Obstructions to homotopy invariance of loop coproduct via parametrised fixed-point theory.
\newblock {\em arXiv preprint arXiv:2407.13662}, 2024.

\bibitem[Nae21]{naef2021string}
Florian Naef.
\newblock {The string coproduct "knows" Reidemeister/Whitehead torsion}.
\newblock {\em arXiv preprint arXiv:2106.11307}, 2021.

\bibitem[NS24]{naef2024simple}
Florian Naef and Pavel Safronov.
\newblock Simple homotopy invariance of the loop coproduct.
\newblock {\em arXiv preprint arXiv:2406.19326}, 2024.

\bibitem[Sul04]{sullivan2004open}
Dennis Sullivan.
\newblock Open and closed string field theory interpreted in classical algebraic topology.
\newblock {\em London Mathematical Society Lecture Note Series}, 308:344, 2004.

\bibitem[Tam10]{tamanoi2010loop}
Hirotaka Tamanoi.
\newblock Loop coproducts in string topology and triviality of higher genus {TQFT} operations.
\newblock {\em Journal of Pure and Applied Algebra}, 214(5):605--615, 2010.

\bibitem[Wah19]{wahl2019invariance}
Nathalie Wahl.
\newblock On the invariance of the string topology coproduct.
\newblock {\em arXiv preprint arXiv:1908.03857}, 2019.

\end{thebibliography}

\end{document}